# MARKOV EQUIVALENCE FOR ANCESTRAL GRAPHS

By R. Ayesha Ali,[1,2] Thomas S. Richardson[3] and Peter Spirtes

*University of Guelph, University of Washington and Carnegie-Mellon University*

Ancestral graphs can encode conditional independence relations that arise in directed acyclic graph (DAG) models with latent and selection variables. However, for any ancestral graph, there may be several other graphs to which it is Markov equivalent. We state and prove conditions under which two maximal ancestral graphs are Markov equivalent to each other, thereby extending analogous results for DAGs given by other authors. These conditions lead to an algorithm for determining Markov equivalence that runs in time that is polynomial in the number of vertices in the graph.

**1. Introduction.** A graphical Markov model is a set of distributions with independence structure described by a graph consisting of vertices and edges. The *independence model* associated with a graph is the set of conditional independence relations encoded by the graph through a *global Markov property*. In general, different graphs may encode the same independence model. In this paper, we consider a particular class of graphs, called ancestral graphs, and characterize when two graphs encode the same sets of conditional independence relations.

The class of ancestral graphs is motivated in the following way. We suppose our observed data were generated by a process represented by a directed acyclic graph (DAG) with a fixed set of variables. The causal interpretation of such a DAG is described by [18] and [14]. However, in general, we may only have observed a subset of these variables in a specific sub-population. Hence, some variables in the underlying DAG are not observed ("latent"),

Received April 2008.
[1]Supported by the National University of Singapore 2004–2005.
[2]Supported in part by NSERC Grant RG 326951-06.
[3]Supported in part by the William and Flora Hewlett Foundation, NSF Grants DMS-99-72008, DMS-05-05865, EPA Grant CR 827970-01-0 and NIH Grant R01 AI032475.

*AMS 2000 subject classifications.* Primary 68T30, 05C75; secondary 68T37.

*Key words and phrases.* Directed acyclic graphs, discriminating path, inducing path, Markov equivalence, polynomial-time algorithm.

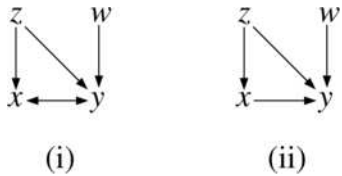

Fig. 1. (i) *A seemingly unrelated regression model and* (ii) *a Markov equivalent DAG model.*

while other variables, specifying the specific sub-population from which our data were sampled, are conditioned upon ("selection variables").

Even though the underlying model is a DAG, the conditional independence structure holding among the observed variables, conditional on the selection variables, cannot always be represented by a DAG containing only the observed variables. For this purpose, the more general class of ancestral graphs is required [see Figure 2(ii) and Definition 2.1]. The statistical models associated with ancestral graphs retain many of the desirable properties that are associated with DAG models.

Like DAGs, two different ancestral graphs can represent the same set of conditional independence relations, and hence distributions. Such graphs are said to be *Markov equivalent*. A graphical characterization of the circumstances under which graphs are Markov equivalent is of importance for several reasons:

- Markov equivalent graphs lead to identical likelihoods because the sets of distributions obeying the Markov property associated with the graphs are the same. Thus, for the purposes of interpreting a model, it is often important to characterize those features that are common to all the graphs in a given class (see [18] and [13]).
- When viewed as a Gaussian path diagram (see [15], Section 8.1), different (maximal) ancestral graphs correspond to different parametrizations of the same Gaussian Markov model. However, some parametrizations may be simpler to fit than others. For example, the model corresponding to the graph in Figure 1(i), in the Gaussian case, is an example of a seemingly unrelated regression (SUR) model (see [23]). In general, there are no closed form expressions for the MLEs for SUR models, iterative fitting methods are required and there may be multiple solutions to the likelihood equations (see [8]). However, the graph in Figure 1(i) is Markov equivalent to Figure 1(ii), which is a DAG. Gaussian DAG models have closed form MLEs, and the likelihood is unimodal (see [12]). Consequently, none of the problems which may arise for general Gaussian SUR models apply to the specific model corresponding to Figure 1(i) (see also [7]).

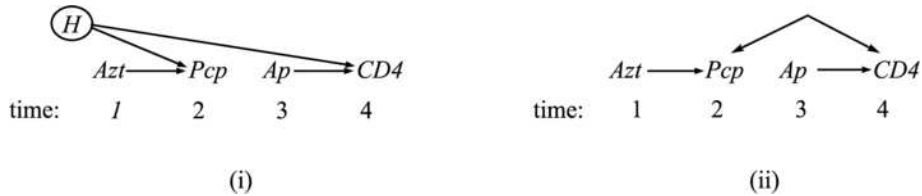

FIG. 2. (i) *A DAG with a latent variable H*. (ii) *The ancestral graph resulting from marginalizing over H includes a bi-directed edge between Pcp and CD4.*

In this paper, we provide necessary and sufficient graphical conditions under which two ancestral graphs are Markov equivalent. Though other characterizations have been given previously in [24] and [19], the criterion given here is the first which leads to an algorithm that runs in time polynomial in the size of the graph. Reference [22] solved the Markov equivalence problem for DAGs. References [2, 3] and [9] solved the problem of representing Markov equivalence classes for DAGs, which we leave for future work.

Section 2 defines the class of ancestral graphs and outlines the motivation for the class. Section 3 contains the main result of the paper. Discussion and relation to prior work are in Section 4. The Appendix contains algorithmic details.

**2. Ancestral graphs.** The basic motivation for ancestral graphs is to enable one to model the independence structure over the observed variables that results from a DAG containing latent or selection variables without explicitly including such variables in the model. To illustrate this, consider the DAG shown in Figure 2(i) in which *Azt*, *Pcp*, *Ap* and *CD4* are observed variables, while $H$ is unobserved. *Azt* and *Ap* represent treatments given to AIDS patients (see Robins [17], Section 2). *Pcp* is an opportunistic infection that often afflicts AIDS patients, and *CD4* can be viewed as a measure of disease progression. Supposing development of *Pcp* was a side-effect of taking *Azt*, then the DAG given in Figure 2(i) incorporates the assumption that *Azt* and *Ap* are both randomized, *Pcp* and *CD4* are responses correlated by underlying health status $H$, and, further, that *Azt* does not affect *CD4*. The DAG implies the following conditional independence relations over the observed variables:

$$Azt \perp\!\!\!\perp Ap, \; CD4, \qquad Ap \perp\!\!\!\perp Azt, \; Pcp.$$

These relations can be derived from the DAG in Figure 2(i) via $d$-separation (see [12] or [22]). Also, note that other valid independence statements, such as $Azt \perp\!\!\!\perp CD4$, can be derived from the two statements given above. The corresponding ancestral graph that represents these same conditional independence relations is shown in Figure 2(ii). (See Section 2.2 for the definition of an ancestral graph and Section 2.3 for the Markov property.) However, there is no DAG on the four observed variables which represents all and only these conditional independence relations.

As this example suggests, bi-directed edges ($\leftrightarrow$) may arise from unobserved parents. Likewise, undirected edges (——) may arise from children that have been conditioned on in the selected sub-population from which the sample is taken (see [4] and [5]). However, bi-directed and undirected edges may also arise in other contexts, where both marginalization and conditioning are present. Reference [16] provides a detailed discussion on the interpretation of edges in an ancestral graph.

2.1. *Basic graphical notation and terminology.* We use the following terminology to describe relations between vertices in a mixed graph $\mathcal{G}$, which may contain three types of edge.

$$\text{If } \left\{\begin{array}{l} a\text{——}b \\ a\leftrightarrow b \\ a\rightarrow b \\ a\leftarrow b \end{array}\right\} \text{ in } \mathcal{G}, \text{ then } a \text{ is a } \left\{\begin{array}{l} neighbor \\ spouse \\ parent \\ child \end{array}\right\} \text{ of } b \quad \text{and} \quad \left\{\begin{array}{l} a \in \text{ne}_{\mathcal{G}}(b) \\ a \in \text{sp}_{\mathcal{G}}(b) \\ a \in \text{pa}_{\mathcal{G}}(b) \\ a \in \text{ch}_{\mathcal{G}}(b) \end{array}\right\}.$$

(For a formal set-theoretic definition of mixed graphs see [15], Appendix.) Two vertices that are connected by some edge are said to be *adjacent*. Note that the three edge types should be considered as distinct symbols, and that all the mixed graphs we consider in this paper are *simple* in that they have at most one edge between each pair of vertices. If there is an edge $a\rightarrow b$ or $a\leftrightarrow b$, then there is said to be *an arrowhead at $b$ on this edge*. Conversely, if there is an edge $a\rightarrow b$ or $a$——$b$, then there is said to be a *tail at $a$*. We also do not allow a vertex to be adjacent to itself. We restrict attention to graphs with finite vertex sets.

A *path* $\boldsymbol{\pi}$ between two vertices $x$ and $y$ in a simple mixed graph $\mathcal{G}$ is a sequence of distinct vertices $\boldsymbol{\pi} = \langle x, v_1, \ldots, v_k, y \rangle$ such that each vertex in the sequence is adjacent to its predecessor and its successor; $x$ and $y$ are the *endpoints* of $\boldsymbol{\pi}$; all other vertices on the path are *nonendpoints* of $\boldsymbol{\pi}$. If $a$ and $b$ are distinct vertices on $\boldsymbol{\pi}$, then the portion of $\boldsymbol{\pi}$ between $a$ and $b$ is called a *section* of $\boldsymbol{\pi}$, denoted $\boldsymbol{\pi}(a,b)$. Note that we use both $\boldsymbol{\pi}(a,b)$ and $\boldsymbol{\pi}(b,a)$ to represent the same section of $\boldsymbol{\pi}$. A path of the form $x\rightarrow\cdots\rightarrow y$, on which every edge is of the form $\rightarrow$, with the arrowheads pointing toward $y$, is a *directed path from $x$ to $y$*. A directed path from $x$ to $y$, together with an edge $y\rightarrow x \in \mathcal{G}$, is called a *directed cycle*.

2.2. *Definition of ancestral graphs.* DAGs are directed graphs in which directed cycles are not permitted. Similarly, certain configurations of edges are not permitted in ancestral graphs:

DEFINITION 2.1. A graph, which may contain undirected (———), directed ($\longrightarrow$) or bi-directed edges ($\leftrightarrow$) is *ancestral* if:

(a) there are no directed cycles;
(b) whenever there is an edge $x \leftrightarrow y$, then there is no directed path from $x$ to $y$, or from $y$ to $x$;
(c) if there is an undirected edge $x$———$y$ then $x$ and $y$ have no spouses or parents.

Conditions (a) and (b) may be summarized by saying that, if $x$ and $y$ are joined by an edge and there is an arrowhead at $x$, then $x$ is *not* an ancestor of $y$; this is the motivation for the term "ancestral."

A vertex $a$ is said to be an *ancestor* of a vertex $b$ if *either* there is a directed path $a \longrightarrow \cdots \longrightarrow b$ from $a$ to $b$ *or* $a = b$. Further, if $a$ is an ancestor of $b$, then $b$ is said to be a *descendant* of $a$.

A vertex $a$ is said to be *anterior* to a vertex $b$ if $a = b$ or there is a path $\mu$ between $a$ and $b$, on which every edge is either of the form $c$———$d$ or $c \longrightarrow d$, with $d$ between $c$ and $b$ on $\mu$; such a path $\mu$ is said to be an *anterior path* from $a$ to $b$. By (c) in Definition 2.1, the configuration $\longrightarrow c$——— never occurs in an ancestral graph; hence, every anterior path takes the form

$$a \text{———} \cdots \text{———} c \longrightarrow \cdots \longrightarrow b,$$

where $a = c$ and $c = b$ are possible. We use $an(x)$, $de(x)$ and $ant(x)$ to denote, respectively, the ancestors of $x$, the descendants of $x$ and the vertices anterior to $x$. We apply these definitions disjunctively to sets. For example,

$$an(X) = \{a \mid a \text{ is an ancestor of } b \text{ for some } b \in X\},$$

$$ant(X) = \{a \mid a \text{ is anterior to } b \text{ for some } b \in X\}.$$

By definition, $X \subseteq an(X) \subseteq ant(X)$. Note that every DAG is an ancestral graph, since clauses (b) and (c) are trivially satisfied.

In the next lemma and elsewhere, we will make use of the shorthand notation $x \text{?}\!\!\longrightarrow y$ to indicate that either $x \longrightarrow y$ or $x \leftrightarrow y$. Similarly, $x \text{?}$———$y$ indicates that either $x \longleftarrow y$ or $x$———$y$, while $x \text{?}$—$\text{?}y$ indicates any edge.

LEMMA 2.2. *Let $a$, $b$, $c$ be vertices in an ancestral graph $\mathcal{G}$ with $a$ and $c$ adjacent. If $a \text{?}\!\!\longrightarrow b \longrightarrow c$, then $a \text{?}\!\!\longrightarrow c$. In particular, if the edge ends at $a$ on the $\langle a, b \rangle$ and $\langle a, c \rangle$ edges differ, then we have $c \longleftarrow a \leftrightarrow b$; otherwise, either $c \longleftarrow a \longrightarrow b$, or $c \leftrightarrow a \leftrightarrow b$.*

We make use of this property in Sections 3.7 and 3.9.

PROOF OF LEMMA 2.2. Suppose, for a contradiction, that there is a tail at $c$ on the $\langle a, c \rangle$ edge. Since, by hypothesis, there is an arrowhead at $b$ on the $\langle b, c \rangle$ edge, $a \text{———} c$ is ruled out by Definition 2.1(c), so $a \prec\!\!\text{—} c$. But then $\mathcal{G}$ violates Definition 2.1(b), since $a\text{?}\!\!\succ b \text{—}\!\!\succ c \text{—}\!\!\succ a$. Hence, $a\text{?}\!\!\succ c$. The conclusion then follows from noting that the configuration $a \text{—}\!\!\succ b \text{—}\!\!\succ c \prec\!\!\!\succ a$ is not ancestral. □

2.3. *The m-separation criterion.* In an ancestral graph, a nonendpoint vertex $v$ on a path is said to be a *collider* if two arrowheads meet at $v$ (i.e., $\text{—}\!\!\succ v \prec\!\!\text{—}$, $\prec\!\!\!\succ v \prec\!\!\!\succ$, $\prec\!\!\!\succ v \prec\!\!\text{—}$ or $\text{—}\!\!\succ v \prec\!\!\!\succ$). All other nonendpoint vertices on a path are *noncolliders* (i.e., $\text{——} v \text{——}$, $\text{——} v \text{—}\!\!\succ$, $\text{—}\!\!\succ v \text{—}\!\!\succ$, $\prec\!\!\text{—} v \text{—}\!\!\succ$, or $\prec\!\!\!\succ v \text{—}\!\!\succ$). These definitions of collider and noncollider are direct extensions of the corresponding definitions for DAGs. A path along which every nonendpoint is a collider is called a *collider path*. A path comprised of 3 vertices is called a *triple*. In an ancestral graph, a triple is either a collider or a noncollider; we refer to this as the *type* of the triple. Hence, if $\langle a, b, c \rangle$ forms a triple, then $\langle c, b, a \rangle$ and $\langle a, b, c \rangle$ are of the same type.

Reference [22] introduced *d-separation*, a set of graphical conditions by which conditional independence relations could be read from a DAG. Reference [15] applied a natural extension of Pearl's $d$-separation criterion, called $m$-separation, to ancestral graphs.

DEFINITION 2.3. Let $a$ and $b$ be distinct vertices in an ancestral graph $\mathcal{G}$, and let $Z$ be a subset of vertices with $a, b \notin Z$. A path $\boldsymbol{\pi}$ between $a$ and $b$ is said to be *m-connecting given $Z$* if the following hold:

(i) no noncollider on $\boldsymbol{\pi}$ is in $Z$; and,
(ii) every collider on $\boldsymbol{\pi}$ is an ancestor of a vertex in $Z$.

Two vertices $a$ and $b$ are said to be *m-separated* given $Z$ in $\mathcal{G}$ if there is no path $m$-connecting $a$ and $b$ given $Z$ in $\mathcal{G}$. Likewise, sets $A$ and $B$ are *m-separated* given $Z$ in $\mathcal{G}$ if, for every pair $a \in A$ and $b \in B$, $a$ and $b$ are $m$-separated given $Z$.

For example, in the ancestral graph in Figure 2(ii), $Azt$ and $Ap$ are $m$-separated given $CD4$. Definition 2.3 is an extension of the original definition of $d$-separation for DAGs in that the notions of "collider" and "noncollider" now allow for bi-directed and undirected edges; the definition of ancestor is unchanged. Furthermore, $d$-separation is equivalent to $m$-separation for DAGs. The following result is useful.

LEMMA 2.4. *In an ancestral graph $\mathcal{G}$, if $\boldsymbol{\pi}$ is a path m-connecting $a$ and $b$ given $Z$, $c$ is on $\boldsymbol{\pi}$ ($a \neq c \neq b$) and there is an arrowhead at $c$ on the section $\boldsymbol{\pi}(a, c)$, then either $c \in \text{an}(Z)$ or $\boldsymbol{\pi}(c, b)$ is a directed path from $c$ to $b$.*

PROOF. Suppose the result is false. Let $c$ be the vertex closest to $b$ satisfying the premise of the lemma but not the conclusion. If $c$ is a collider on $\pi$, then, by definition of $m$-connection, $c \in \text{an}(Z)$ which is a contradiction. Let $d$ be the vertex after $c$ on $\pi(c,b)$. If $c$ is a noncollider on $\pi$ then, by Definition 2.1(c), $c \longrightarrow d$. If $d \in \text{an}(Z)$ or $\pi(d,b)$ forms a directed path from $d$ to $b$, then, clearly, $c$ satisfies the conclusion of the lemma. But, if $d \notin \text{an}(Z)$ and $\pi(d,b)$ is not a directed path to $b$, then $d$ satisfies the premise of the lemma (and hence, $c$ is not the closest such vertex to $b$), again a contradiction. $\square$

2.4. *Formal independence models.* An *independence model* over a finite set $V$ is a set $\mathfrak{I}$ of ternary relations $\langle X, Y \mid Z \rangle$ where $X$, $Y$ and $Z$ are disjoint subsets of $V$, while $X$ and $Y$ are not empty; the first two arguments are treated symmetrically, so that $\langle X, Y \mid Z \rangle \in \mathfrak{I}$ iff $\langle Y, X \mid Z \rangle \in \mathfrak{I}$. The interpretation of $\langle X, Y \mid Z \rangle \in \mathfrak{I}$ is that $X$ and $Y$ are independent given $Z$ [see [20], Chapter 2]. The independence model associated with an ancestral graph, $\mathfrak{I}_m(\mathcal{G})$, is defined via $m$-separation as follows

$$\mathfrak{I}_m(\mathcal{G}) \equiv \{\langle X, Y \mid Z \rangle | X \text{ is } m\text{-separated from } Y \text{ given } Z \text{ in } \mathcal{G}\}.$$

The independence relations in $\mathfrak{I}_m(\mathcal{G})$ comprise the *global Markov property* for $\mathcal{G}$.

2.5. *Probability distributions obeying a formal independence model.* We associate a set of probability distributions with a formal independence model $\mathfrak{I}$ by using the finite set $V$ to index a collection of random variables $(X_\nu)_{\nu \in V}$ taking values in probability spaces $(\Omega_\nu)_{\nu \in V}$. In all the examples we consider, the probability spaces are either real finite-dimensional vector spaces or finite discrete sets. For $A \subseteq V$, we let $\Omega_A \equiv \times_{v \in A} (\Omega_\nu)$, $\Omega \equiv \Omega_V$ and $X_A \equiv (X_\nu)_{\nu \in A}$. We will assume the existence of regular conditional probability measures throughout.

A distribution $P$ on $\Omega$ is said to *obey the independence model* $\mathfrak{I}$ over $V$ if, for all disjoint sets $A, B, Z$ ($A$ and $B$ are not empty),

$$\langle A, B \mid Z \rangle \in \mathfrak{I} \quad \Longrightarrow \quad A \perp\!\!\!\perp B \mid Z[P],$$

where we have used the ($\perp\!\!\!\perp$) notation of [6], and the usual shorthand that $A$ denotes both a vertex set and the random variable $X_A$. Thus, a distribution $P$ obeys $\mathfrak{I}_m(\mathcal{G})$ if, for all disjoint subsets of $V$, say $X$, $Y$, $Z$ ($X$ and $Y$ not empty),

$$X \text{ is } m\text{-separated from } Y \text{ given } Z \text{ in } \mathcal{G} \quad \Longrightarrow \quad X \perp\!\!\!\perp Y \mid Z[P].$$

Note that, if $P$ obeys $\mathfrak{I}$, there still may be independence relations that are not in $\mathfrak{I}$ that also hold in $P$.

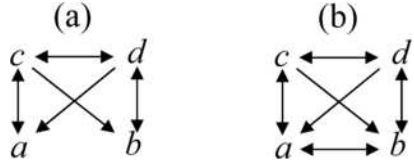

Fig. 3. (a) *The path $\langle a, c, d, b \rangle$ is an example of an inducing path in an ancestral graph.* (b) *A maximal ancestral graph Markov equivalent to* (a).

2.6. *Marginalizing and conditioning.* In Section 4.1 of [15] operations of marginalizing and conditioning are introduced for formal independence models. If $P$ obeys $\mathfrak{I}$ and $\mathfrak{I}^*$ is the independence model obtained by formally marginalizing over variables in $L$ and conditioning on variables in $S$, then $P(X_{V \setminus (L \cup S)} \mid X_S)$ obeys the independence model $\mathfrak{I}^*$ [$P(X_S)$ a.e.] (see Theorem 7.1 of [15], Appendices A and B of [10]).

In Section 4.2 of [15], a graphical transformation corresponding to marginalizing and conditioning is given such that the independence model associated with the transformed graph is the independence model obtained by marginalizing and conditioning the independence model $\mathfrak{I}_m(\mathcal{G})$ of the original graph (see Theorem 4.18 in [15]). Thus, in particular, if $\mathcal{G}$ is a DAG with observed variables $O$, latent variables $L$ and selection variables $S$, then the ancestral graph formed by the graphical transformation applied to $\mathcal{G}$ represents those conditional independence relations implied to hold among the observed variables $O$, conditional on the selection variables [$P(X_S)$ a.e.].

**3. Markov equivalence.** We introduce the following.

DEFINITION 3.1. Two ancestral graphs $\mathcal{G}_1$ and $\mathcal{G}_2$ with the same vertex set are said to be *Markov equivalent*, denoted $\mathcal{G}_1 \sim \mathcal{G}_2$, if for all disjoint sets $A$, $B$, $Z$ ($A$, $B$ not empty), $A$ and $B$ are $m$-separated given $Z$ in $\mathcal{G}_1$ if and only if $A$ and $B$ are $m$-separated given $Z$ in $\mathcal{G}_2$; that is, $\mathfrak{I}_m(\mathcal{G}_1) = \mathfrak{I}_m(\mathcal{G}_2)$.

The graphs in Figure 3 are Markov equivalent, as are $\mathcal{G}_1$ and $\mathcal{G}_2$ in Figure 4. The set of all ancestral graphs that encode the same set of conditional independence statements forms a *Markov equivalence class*.

3.1. *Markov equivalence for DAGs.* References [9] and [22] gave simple graphical conditions for determining whether two DAGs are Markov equivalent. A triple of vertices $\langle a, b, c \rangle$ is said to be *unshielded* if $a$ and $c$ are not adjacent and *shielded* otherwise. (A triple is defined in Section 2.3.)

THEOREM 3.2. *Two DAGs are Markov equivalent if and only if they have the same adjacencies and the same unshielded colliders.*

That two Markov equivalent DAGs have the same adjacencies is a direct consequence of the fact that DAGs satisfy a pairwise Markov property.

PROPOSITION 3.3 ([12], page 50). *In a DAG $\mathcal{D}$, if $a$ and $b$ are not adjacent and $b \notin \mathrm{an}(a)$, then $a$ is d-separated from $b$ by $V \setminus (\mathrm{de}(b) \cup \{a\})$.*

This is a consequence of the local Markov property for DAGs [12], applied to $b$, which implies that $b$ is d-separated from $V \setminus (\mathrm{pa}(b) \cup \mathrm{de}(b))$ by $\mathrm{pa}(b)$; $b \notin \mathrm{an}(a)$ implies $a \notin \mathrm{de}(b)$ and $\mathrm{pa}(b) \subseteq V \setminus (\mathrm{de}(b) \cup \{a\})$. Note that, by acyclicity, for any pair $a$, $b$, either $b \notin \mathrm{an}(a)$ or $a \notin \mathrm{an}(b)$. Consequently, in a DAG, every missing edge implies a conditional independence between the nonadjacent vertices. In general, no such pairwise property holds for ancestral graphs. For example, there is no set that $m$-separates $a$ and $b$ in the graph in Figure 3(a). This motivates the following section.

### 3.2. Maximal ancestral graphs.

DEFINITION 3.4. An ancestral graph $\mathcal{G}$ is said to be *maximal* if, for every pair of nonadjacent vertices $(a, b)$, there exists a set $Z$ $(a, b \notin Z)$ such that $a$ and $b$ are $m$-separated conditional on $Z$.

These graphs are maximal in the sense that no additional edge may be added to the graph without changing the associated independence model. In a nonmaximal ancestral graph two nonadjacent vertices $a$ and $b$, for which no $m$-separating set $Z$ exists, will be joined by an *inducing path*.

DEFINITION 3.5. An *inducing path* $\boldsymbol{\pi}$ between vertices $a$ and $b$ in an ancestral graph $\mathcal{G}$ is a path on which every nonendpoint vertex is both a collider on $\boldsymbol{\pi}$ and an ancestor of at least one of the endpoints, $a$, $b$.

For a proof, see [15], Corollary 4.3, where the definition given here is termed a "primitive" inducing path; the concept was introduced by Verma and Pearl [21]. Note that, strictly speaking, an inducing "path" $\boldsymbol{\pi} = \langle a, v_1, \ldots, v_k, b \rangle$ is a collection of paths: the collider path $\boldsymbol{\pi}$, together with directed paths from each vertex $v_i$, $1 \leq i \leq k$, to one of the endpoints. The name "inducing" path refers to the fact that given any set $Z$ $(a,b \notin Z)$ $\boldsymbol{\pi}$ is $m$-connecting given $Z$. If there is some vertex $v_i \notin \mathrm{an}(Z)$, then there is an $m$-connecting path involving one or more of the directed paths, otherwise the path $\boldsymbol{\pi}$ itself is $m$-connecting.

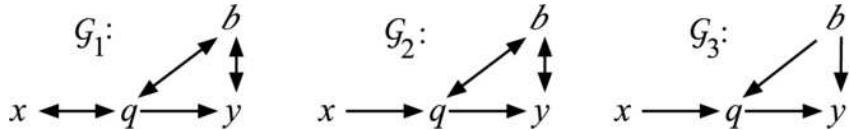

FIG. 4. $\mathcal{G}_1$, $\mathcal{G}_2$, $\mathcal{G}_3$ *have the same adjacencies and the same unshielded colliders, but $\mathcal{G}_1$ and $\mathcal{G}_3$ are not Markov equivalent.* $\boldsymbol{\pi} = \langle x, q, b, y \rangle$ *forms a discriminating path for $b$ in every graph.*

Figure 3(a) shows an example of a nonmaximal ancestral graph. The path $\langle a, c, d, b \rangle$ forms an inducing path between $a$ and $b$. By adding the bi-directed edge $a \leftrightarrow b$, the graph is made maximal without changing the associated independence model (which is empty), as shown in Figure 3(b). As is the case in this example, in general, if $\boldsymbol{\pi} = \langle a, v_1, \ldots, v_k, b \rangle$ is an inducing path, then only a bi-directed edge $a \leftrightarrow b$ may be added while obeying Definition 2.1. Since there are arrowheads present at $a$ and $b$, adding an undirected edge is ruled out by (c); adding a directed edge would violate (b) since we would either have $a \leftrightarrow v_1 \rightarrow \cdots \rightarrow b \rightarrow a$ or $b \leftrightarrow v_k \rightarrow \cdots \rightarrow a \rightarrow b$.

By [15], Theorem 5.1, for every nonmaximal ancestral graph $\mathcal{G}$ there is a unique maximal ancestral graph $\bar{\mathcal{G}}$ of which it is a subgraph; in fact, $\bar{\mathcal{G}} = \mathcal{G}[_\varnothing^\varnothing$ and thus $\bar{\mathcal{G}}$ may be constructed in polynomial time. Consequently, the problem of characterizing Markov equivalence for ancestral graphs naturally reduces to that of characterizing equivalence in the case where both graphs are maximal. Except where noted, in the remainder of this paper, we will restrict attention to maximal ancestral graphs (MAGs).

3.3. *Necessary conditions for Markov equivalence.*

PROPOSITION 3.6. *If $\mathcal{G}_1$, $\mathcal{G}_2$ are MAGs and $\mathcal{G}_1 \sim \mathcal{G}_2$, then $\mathcal{G}_1$ and $\mathcal{G}_2$ have the same adjacencies and unshielded colliders.*

PROOF. Since $\mathcal{G}_1$ is maximal, for each pair of nonadjacent vertices $(x, y)$ in $\mathcal{G}_1$, there is some set $Z$ such that $x$ and $y$ are $m$-separated given $Z$ in $\mathcal{G}_1$. If $x$ and $y$ are adjacent in $\mathcal{G}_2$, then they are not $m$-separated by $Z$, contradicting $\mathcal{G}_1 \sim \mathcal{G}_2$. So, adjacencies in $\mathcal{G}_1$ are a subset of those in $\mathcal{G}_2$. By a symmetric argument, the adjacencies in $\mathcal{G}_2$ are a subset of those in $\mathcal{G}_1$.

Suppose, for a contradiction, that $\langle a, b, c \rangle$ is an unshielded collider in $\mathcal{G}_1$ but not in $\mathcal{G}_2$. Since $\mathcal{G}_1$ is maximal, for some set $Z$, $a$ and $c$ are $m$-separated by $Z$, and $b \notin Z$. If $\langle a, b, c \rangle$ is a noncollider in $\mathcal{G}_2$ then $a$ and $c$ are $m$-connected given $Z$, which is a contradiction. Hence, every unshielded collider in $\mathcal{G}_1$ is present in $\mathcal{G}_2$. The conclusion follows by symmetry. □

An important consequence of this proposition is that if $\mathcal{G}_1$ and $\mathcal{G}_2$ are maximal and Markov equivalent, then a sequence of vertices forming a path in $\mathcal{G}_1$ also forms a path in $\mathcal{G}_2$ and vice-versa, though the edge-types on these paths may differ. Consequently, when $\mathcal{G}_1 \sim \mathcal{G}_2$, we will often refer to the path $\boldsymbol{\pi}^*$ in $\mathcal{G}_2$ *corresponding to* a given path $\boldsymbol{\pi}$ in $\mathcal{G}_1$.

A key difference between DAGs and MAGs is that having the same adjacencies and the same unshielded colliders, though necessary, are no longer sufficient for Markov equivalence. Consider the graphs shown in Figure 4. $\mathcal{G}_1$ and $\mathcal{G}_3$ contain the same adjacencies and the same unshielded colliders, but these two graphs are not Markov equivalent to each other. In $\mathcal{G}_1$, $x$ is $m$-separated from $y$ given $q$; but according to $\mathcal{G}_3$, $x$ is $m$-connected to $y$ given $q$. In fact, in any graph Markov equivalent to $\mathcal{G}_1$, $\langle q, b, y \rangle$ forms a shielded collider. (There is only one such graph, $\mathcal{G}_2$, so $\{\mathcal{G}_1, \mathcal{G}_2\}$ forms a Markov equivalence class.) However, in general, it is clearly not necessary that two graphs have *all* of the same shielded colliders in order for them to be Markov equivalent. Much of the remainder of this paper will focus on identifying the "relevant" set of colliders for judging Markov equivalence. The main result of this paper follows.

THEOREM 3.7. *If $\mathcal{G}_1$, $\mathcal{G}_2$ are MAGs, then $\mathcal{G}_1 \sim \mathcal{G}_2$ if and only if $\mathcal{G}_1$ and $\mathcal{G}_2$ have the same adjacencies and the same colliders with order.*

The set of "colliders with order" within a graph is defined recursively in Definition 3.11 in the next section. The proof concludes in Section 3.10.

3.4. *Discriminating paths in maximal ancestral graphs.* A *discriminating path*, if present in two Markov equivalent MAGs, implies that a certain shielded triple will be of the same type in both graphs.

DEFINITION 3.8 [18]. A path $\boldsymbol{\pi} = \langle x, q_1, \ldots, q_p, b, y \rangle$ ($p \geq 1$) is a *discriminating path for* $\langle q_p, b, y \rangle$ in a MAG $\mathcal{G}$ if:

(i) $x$ is not adjacent to $y$, and,
(ii) every vertex $q_i$ ($1 \leq i \leq p$) is a collider on $\boldsymbol{\pi}$, and a parent of $y$.

We will often refer to a section $\boldsymbol{\pi}(x, y)$ of some path $\boldsymbol{\pi}$ as *a discriminating path for $b$*, thereby implicitly specifying the triple $\langle q_p, b, y \rangle = \boldsymbol{\pi}(q_p, y)$. By convention, we order the endpoints of the discriminating path so it is the second endpoint (in this case, $y$) which is in the discriminated triple. We are free to order $x$ and $y$ in this way, since, in our notation, $\boldsymbol{\pi}(x, y)$ and $\boldsymbol{\pi}(y, x)$ represent the same section of $\boldsymbol{\pi}$ (see page 4).

The paths $\langle x, q, b, y \rangle$ in $\mathcal{G}_1$, $\mathcal{G}_2$ and $\mathcal{G}_3$ from Figure 4 are examples of discriminating paths for $b$. Like an inducing path, a discriminating "path" $\boldsymbol{\pi} = \langle x, q_1, \ldots, q_p, b, y \rangle$ is, in fact, a collection of paths

$$x \mathbin{?\!\!\rightarrow} q_1 \leftrightarrow \cdots \leftrightarrow q_j \rightarrow y \qquad (1 \leq j \leq p),$$

$$x \mathbin{?\!\!\rightarrow} q_1 \leftrightarrow \cdots \leftrightarrow q_p \leftarrow\mathbin{?} b \mathbin{?\!\!\rightarrow} y,$$

together with the (additional) requirement that the endpoints $x$ and $y$ are not adjacent. Consider a discriminating path $\boldsymbol{\pi} = \langle x, q_1, \ldots, q_p, b, y \rangle$ in an ancestral graph $\mathcal{G}$. If a given set $Z$ $(x, y \notin Z)$ does not contain all vertices $q_i, 1 \leq i \leq p$, then, for some $j$, $q_j \notin Z$ and for all $k < j$, $q_k \in Z$, so that the path $\langle x, q_1, \ldots, q_j, y \rangle$ $m$-connects $x$ and $y$ given $Z$ (because $q_1, \ldots, q_{j-1}$ are colliders and $q_j$ is a noncollider); see Figure 5. Hence, if $Z$ $m$-separates $x$ and $y$ then $\{q_1, \ldots, q_p\} \subseteq Z$. Consequently, if $b$ is a collider on the path $\boldsymbol{\pi}$ in the graph $\mathcal{G}$ and $Z$ $m$-separates $x$ and $y$, then $b \notin Z$; otherwise, the path $\boldsymbol{\pi}$ would $m$-connect $x$ and $y$, since every nonendpoint vertex on $\boldsymbol{\pi}$ would be a collider and in $Z$. Conversely, if $b$ is a noncollider on the path $\boldsymbol{\pi}$, then $b$ is a member of any set $Z$ that $m$-separates $x$ and $y$.

Thus, whenever $\langle x, q_1, \ldots, q_p, b, y \rangle$ forms a discriminating path in $\mathcal{G}$, then $b$ is a collider [noncollider] if and only if every set $Z$ $m$-separating $x$ and $y$ is such that $b \notin Z$ [$b \in Z$]. It follows that if $\mathcal{G}^* \sim \mathcal{G}$ and the path corresponding to $\boldsymbol{\pi}$, say $\boldsymbol{\pi}^*$, also forms a discriminating path for $b$ in $\mathcal{G}^*$, then $b$ is a collider on $\boldsymbol{\pi}^*$ (in $\mathcal{G}^*$) if and only if $b$ is a collider on $\boldsymbol{\pi}$ (in $\mathcal{G}$). Thus, we have proved the following.

LEMMA 3.9. *Let $\boldsymbol{\pi} = \langle x, q_1, \ldots, q_p, b, y \rangle$ be a discriminating path for $b$ in the MAG $\mathcal{G}$. If $\mathcal{G}^*$ is a MAG, $\mathcal{G}^* \sim \mathcal{G}$, and the corresponding path $\boldsymbol{\pi}^*$ forms a discriminating path for $b$ in $\mathcal{G}^*$, then $b$ is a collider on $\boldsymbol{\pi}$ in $\mathcal{G}$ if and only if $b$ is a collider on $\boldsymbol{\pi}^*$ in $\mathcal{G}^*$.*

Thus, in general, even though $q_p$ and $y$ are adjacent, $\langle q_p, b, y \rangle$ is "discriminated" by the path $\boldsymbol{\pi}$ to be of the same type (collider or noncollider) on the corresponding path in any graph $\mathcal{G}^*$ Markov equivalent to $\mathcal{G}$ in which *the corresponding path $\boldsymbol{\pi}^*$ also forms a discriminating path*. Though discriminating paths can exist in DAGs, they are not important for determining Markov equivalence, because such paths always discriminate noncolliders (see $\mathcal{G}_3$ in Figure 4). If $\langle x, q, b \rangle$ forms a collider, then since there are no bi-directed edges in a DAG, it follows that $b$ is a parent of $q$.

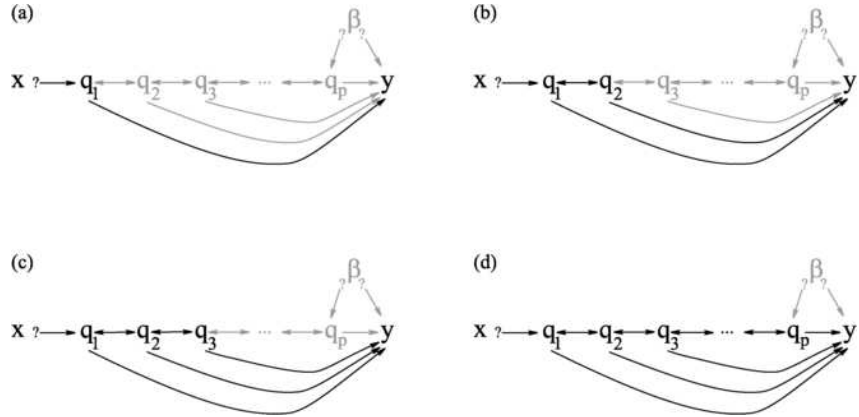

FIG. 5. *The unshielded noncolliders $\langle x, q_1, y \rangle$ and the sequence of discriminating paths for the noncolliders $\langle q_{j-1}, q_j, y \rangle$ ($1 < j \leq p$). See Lemma 3.10 and Corollary 3.14.*

The following lemma gives a sufficient condition under which the path $\pi^*$ corresponding to a discriminating path $\pi$ in a MAG $\mathcal{G}$ will also be discriminating in another Markov equivalent MAG $\mathcal{G}^*$.

LEMMA 3.10. *If $\pi = \langle x, q_1, \ldots, q_p, b, y \rangle$ is a discriminating path in a MAG $\mathcal{G}$, then, in any MAG $\mathcal{G}^*$ with $\mathcal{G}^* \sim \mathcal{G}$ in which the $q_i$ are colliders on the corresponding path $\pi^*$, the edges between $q_i$ and $y$ in $\mathcal{G}^*$ are of the form $q_i \rightarrowtail y$, $(1 \leq i \leq p)$.*

PROOF. The proof proceeds by induction on $i$. First, consider the $\langle q_1, y \rangle$ edge in $\mathcal{G}^*$. If there is an arrowhead at $q_1$, then $\langle x, q_1, y \rangle$ forms an unshielded collider in $\mathcal{G}^*$ but an unshielded noncollider in $\mathcal{G}$. But then, by Proposition 3.6, $\mathcal{G}$ and $\mathcal{G}^*$ are not Markov equivalent, which is a contradiction. Since $x \mathbin{?\!\rightarrowtail} q_1 \mathbin{\text{---}?} y$, but $q_1 \text{---} y$ is ruled out by Definition 2.1(c), we have $q_1 \rightarrowtail y$ in $\mathcal{G}^*$.

Suppose that $q_j \rightarrowtail y$ for $1 \leq j < i$ in $\mathcal{G}^*$. Then, the path $\langle x, q_1, \ldots, q_i, y \rangle$, $i \leq p$ forms a discriminating path for $q_i$ in both $\mathcal{G}^*$ and $\mathcal{G}$. If $q_i \mathbin{\leftarrowtail\!?} y$ in $\mathcal{G}^*$, then $\langle q_{i-1}, q_i, y \rangle$ forms a collider in $\mathcal{G}^*$ but a noncollider in $\mathcal{G}$. But then, by Lemma 3.9, we have $\mathcal{G} \not\sim \mathcal{G}^*$, which is a contradiction. Since $q_{i-1} \leftrightarrow q_i \text{---}? y$, but $q_i \text{---} y$ is ruled out by Definition 2.1(c), we have $q_i \rightarrowtail y$ in $\mathcal{G}^*$ as required. □

One might hope that, if $\mathcal{G}_1 \sim \mathcal{G}_2$, then $\mathcal{G}_1$ and $\mathcal{G}_2$ would have the same discriminating paths. Unfortunately, this is not the case. It is possible for a path $\pi$ to be discriminating in $\mathcal{G}$, and yet the corresponding path $\pi^*$ not be discriminating in $\mathcal{G}^*$ even though $\mathcal{G} \sim \mathcal{G}^*$. Hence, the premise in Lemma 3.9 will not hold for all pairs of Markov equivalent graphs. Thus, the fact that a noncollider is discriminated by a path in $\mathcal{G}$ does not mean that it will be present in every graph Markov equivalent to $\mathcal{G}$.

Consider the example given by the two graphs in Figure 6(i). Note that $q$ is a collider on the path $\langle x, q, b, y \rangle$ in $\mathcal{G}_1$, but not in $\mathcal{G}_2$; $\langle x, q, b, y \rangle$ forms a discriminating path in $\mathcal{G}_1$, but not in $\mathcal{G}_2$, though $\mathcal{G}_1 \sim \mathcal{G}_2$. Hence, although $\langle q, b, y \rangle$ is a noncollider in any graph Markov equivalent to $\mathcal{G}_1$ in which $\langle x, q, b, y \rangle$ forms a discriminating path for $b$, $\langle q, b, y \rangle$ need not be a noncollider in graphs such as $\mathcal{G}_2$, where the corresponding path is not discriminating for $b$.

However, we conjecture that if a collider is discriminated by some path in $\mathcal{G}$, then this collider will be present in every graph $\mathcal{G}^*$ Markov equivalent to $\mathcal{G}$, regardless of whether there is a discriminating path for this collider in $\mathcal{G}^*$ or not. For example, the collider $\langle q, b, y \rangle$ in the graph $\mathcal{G}_1$, shown in Figure 6(ii), is present in every graph Markov equivalent to $\mathcal{G}_1$, even though the path $\langle x, q, b, y \rangle$ does not always form a discriminating path, as in $\mathcal{G}_2$, shown in Figure 6(ii).

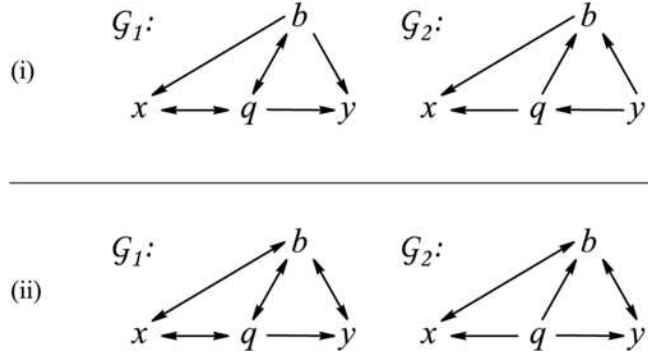

FIG. 6. *Two examples of maximal ancestral graphs that are Markov equivalent where $\langle x, q, b, y \rangle$ forms a discriminating path in $\mathcal{G}_1$, but not in $\mathcal{G}_2$.*

The results in this section present a dilemma; it is clear that discriminating paths, when present in both graphs, lead directly to necessary conditions for Markov equivalence. However, a discriminating path for a given triple may not be present in all graphs within a Markov equivalence class. We avoid this problem by identifying, via a recursive definition, a sub-class of discriminating paths and associated triples (those "with order") that are always present, and by showing that, in conjunction with the conditions in Proposition 3.6, these triples provide sufficient conditions for determining Markov equivalence.

DEFINITION 3.11. Let $\mathfrak{D}_i$ ($i \geq 0$) be the set of *triples of order $i$* in a MAG $\mathcal{G}$, defined recursively as follows:

*Order 0.* A triple $\langle a, b, c \rangle \in \mathfrak{D}_0$ if $a$ and $c$ are not adjacent.
*Order $i+1$.* A triple $\langle a, b, c \rangle \in \mathfrak{D}_{i+1}$ if
    (1) for all $j < i+1$, $\langle a, b, c \rangle \notin \mathfrak{D}_j$, and,
    (2) there is a discriminating path $\langle x, q_1, \ldots, q_p, b, y \rangle$ for $b$ with either $\langle a, b, c \rangle = \langle q_p, b, y \rangle$ or $\langle a, b, c \rangle = \langle y, b, q_p \rangle$ and the $p$ colliders
$$\langle x, q_1, q_2 \rangle, \ldots, \langle q_{p-1}, q_p, b \rangle \in \bigcup_{j \leq i} \mathfrak{D}_j.$$

If $\langle a, b, c \rangle \in \mathfrak{D}_i$ then the triple is said to *have order $i$*. If a triple has order $i$ for some $i$, then we will say that the triple *has order*. A *discriminating path is said to have order $i$* if, excepting $\langle q_p, b, y \rangle$, every collider on the path has order at most $i - 1$, and at least one collider has order $i - 1$.

For example, in every graph in Figure 4, the triple $\langle x, q, b \rangle$ has order 0, while $\langle q, b, y \rangle$ has order 1. It is important to note that not every triple in a graph will have an order. For example, in all the graphs in Figure 6, the triples $\langle x, q, b \rangle$ and $\langle q, b, y \rangle$ do not have order. However, it is possible for a triple without order to be of the same type (collider or noncollider) in every graph in the Markov equivalence class, such as triple $\langle q, b, y \rangle$ in Figure 6(ii). Note that the order (if any) of a shielded triple is the minimum of the orders of all discriminating paths (with order) for that triple.

We now show that a necessary condition for two graphs to be Markov equivalent is that they have the same colliders with order.

PROPOSITION 3.12. *If $\langle a, b, c \rangle$ has order $r$ in a MAG $\mathcal{G}$, then $\langle a, b, c \rangle$ has order $r$ in any MAG $\mathcal{G}^*$, with $\mathcal{G}^* \sim \mathcal{G}$, and, further, $\langle a, b, c \rangle$ is a collider in $\mathcal{G}$ if and only if $\langle a, b, c \rangle$ is a collider in $\mathcal{G}^*$.*

PROOF. The proof is by induction on $r$, the order of $\langle a, b, c \rangle$. For $r = 0$, the result follows from Proposition 3.6. For $r > 0$, by Definition 3.11, there exists a discriminating path $\boldsymbol{\pi} = \langle q_0, \ldots, q_p = a, b, c \rangle$ or $\langle q_0, \ldots, q_p = c, b, a \rangle$ in $\mathcal{G}$ such that, with the possible exception of $\langle a, b, c \rangle$, every other triple on $\boldsymbol{\pi}$ is a collider and has order less than $r$. By the induction hypothesis, in $\mathcal{G}^*$ these triples have the same order as in $\mathcal{G}$ and also form colliders. By Lemma 3.10, since the $q_i$'s ($i > 0$) are colliders on the corresponding path $\boldsymbol{\pi}^*$ in $\mathcal{G}^*$, $q_i \rightarrowtail y$ ($1 \leq i \leq p$) in $\mathcal{G}^*$. Thus, $\boldsymbol{\pi}^*$ also forms a discriminating path in $\mathcal{G}^*$, and so $\langle a, b, c \rangle$ has order at most $r$ in $\mathcal{G}^*$. However, if $\langle a, b, c \rangle$ has order less than $r$ in $\mathcal{G}^*$, then, by the inductive hypothesis (applied to $\mathcal{G}^*$), $\langle a, b, c \rangle$ will have lower order than $r$ in $\mathcal{G}$, contrary to assumption. Thus, $\langle a, b, c \rangle$ has order $r$ in $\mathcal{G}^*$. The result follows by Lemma 3.9. □

LEMMA 3.13. *If MAGs $\mathcal{G}_1$ and $\mathcal{G}_2$ have the same adjacencies and are such that:*

   (i) *every collider with order in $\mathcal{G}_1$ is a collider in $\mathcal{G}_2$, and*
   (ii) *every collider with order in $\mathcal{G}_2$ is a collider in $\mathcal{G}_1$,*

*then, for all $r \geq 0$, $\langle a, b, c \rangle$ is a collider [noncollider] with order $r$ in $\mathcal{G}_1$ iff $\langle a, b, c \rangle$ is a collider [noncollider] with order $r$ in $\mathcal{G}_2$.*

It will follow from Lemma 3.13 and Theorem 3.7 below that conditions (i) and (ii), together with the same adjacencies, are sufficient for Markov equivalence.

PROOF OF LEMMA 3.13. We argue by induction for each order $r$.

($r = 0$). Suppose $\langle a, b, c \rangle$ is a triple of order 0 in $\mathcal{G}_1$ [$\mathcal{G}_2$]. By Definition 3.11, $a$ and $c$ are not adjacent in $\mathcal{G}_1$ [$\mathcal{G}_2$]. Hence, $a$ and $c$ are not adjacent in $\mathcal{G}_2$ [$\mathcal{G}_1$], so $\langle a,b,c \rangle$ has order 0 in $\mathcal{G}_2$ [$\mathcal{G}_1$]. If $\langle a,b,c \rangle$ forms a collider in $\mathcal{G}_1$ [$\mathcal{G}_2$], then, by (i) [(ii)], it forms a collider (with order 0) in $\mathcal{G}_2$ [$\mathcal{G}_1$]. Conversely, if $\langle a,b,c \rangle$ forms a noncollider in $\mathcal{G}_1$ [$\mathcal{G}_2$], then, since it also has order 0 in $\mathcal{G}_2$ [$\mathcal{G}_1$], by (ii) [(i)], $\langle a,b,c \rangle$ cannot be a collider in $\mathcal{G}_2$ [$\mathcal{G}_1$].

($r > 0$). Suppose the result holds for all $s < r$. If $\langle a,b,c \rangle$ is a triple with order $r$ in $\mathcal{G}_1$ [$\mathcal{G}_2$], then there is a discriminating path $\boldsymbol{\mu} = \langle q_0, q_1, \ldots, q_p, b, y \rangle$, where either $q_p = a$ and $y = c$, or $q_p = c$ and $y = a$, and each collider $q_i$ ($1 \leq i \leq p$) on $\boldsymbol{\mu}$ has order less than $r$ by Definition 3.11. By the induction hypothesis, each collider $q_i$ is also a collider on the corresponding path $\boldsymbol{\mu}^*$ in $\mathcal{G}_2$ [$\mathcal{G}_1$] with the same order as in $\mathcal{G}_1$ [$\mathcal{G}_2$].

We claim that $\boldsymbol{\mu}^*$ also forms a discriminating path in $\mathcal{G}_2$ [$\mathcal{G}_1$]. Since we have $q_0 ?\!\!\rightarrow q_1 \leftrightarrow \cdots \leftrightarrow q_p \leftarrow\!\!? b$ in $\mathcal{G}_2$ [$\mathcal{G}_1$], it suffices to show that $q_j \rightarrow y$ ($1 \leq j \leq p$) in $\mathcal{G}_2$ [$\mathcal{G}_1$] (see Figure 5). Triple $\langle q_0, q_1, y \rangle$ is a noncollider with order 0 in $\mathcal{G}_1$ [$\mathcal{G}_2$] because $q_0$ and $y$ are not adjacent. Hence, by the inductive hypothesis, $\langle q_0, q_1, y \rangle$ is a noncollider (with order 0) in $\mathcal{G}_2$ [$\mathcal{G}_1$]. Further, by Definition 2.1(c), $q_1 \rightarrow y$ in $\mathcal{G}_2$ [$\mathcal{G}_1$], because $q_0 ?\!\!\rightarrow q_1$. Arguing inductively, assume that $q_i \rightarrow y$ ($1 \leq i < j$) in $\mathcal{G}_2$ [$\mathcal{G}_1$] so that $\langle q_0, q_1, \ldots, q_j, y \rangle$ forms a discriminating path with order at most $r$ for $\langle q_{j-1}, q_j, y \rangle$ in both graphs. Consequently, if $\langle q_{j-1}, q_j, y \rangle$ formed a collider in $\mathcal{G}_2$ [$\mathcal{G}_1$], then $\langle q_{j-1}, q_j, y \rangle$ would be a collider with order at most $r$ in $\mathcal{G}_2$ [$\mathcal{G}_1$] but a noncollider in $\mathcal{G}_1$ [$\mathcal{G}_2$], contrary to (ii) [(i)]. Since $q_{j-1} ?\!\!\rightarrow q_j$ and $\langle q_{j-1}, q_j, y \rangle$ forms a noncollider, by Definition 2.1(c), $q_j \rightarrow y$ in $\mathcal{G}_2$ [$\mathcal{G}_1$].

Hence, $\boldsymbol{\mu}^*$ forms a discriminating path with order at most $r$ in $\mathcal{G}_2$ [$\mathcal{G}_1$], so $\langle a,b,c \rangle$ has order at most $r$ in $\mathcal{G}_2$ [$\mathcal{G}_1$]. However, if $\langle a,b,c \rangle$ has order less than $r$ in $\mathcal{G}_2$ [$\mathcal{G}_1$], then, by the inductive hypothesis, $\langle a,b,c \rangle$ will have lower order than $r$ in $\mathcal{G}_1$ [$\mathcal{G}_2$], contrary to assumption. Thus, $\langle a,b,c \rangle$ has order $r$ in both graphs.

Now, if $\langle a,b,c \rangle$ is a collider in $\mathcal{G}_1$ [$\mathcal{G}_2$], then, by (i) [(ii)], $\langle a,b,c \rangle$ is also a collider in $\mathcal{G}_2$ [$\mathcal{G}_1$]. Conversely, if $\langle a,b,c \rangle$ is a noncollider in $\mathcal{G}_1$ [$\mathcal{G}_2$], then it cannot be a collider in $\mathcal{G}_2$ [$\mathcal{G}_1$] as that would violate (ii) [(i)]. □

COROLLARY 3.14. *If MAGs $\mathcal{G}_1$ and $\mathcal{G}_2$ have the same adjacencies and $\langle a,b,c \rangle$ is a collider with order in $\mathcal{G}_1$ iff $\langle a,b,c \rangle$ is a collider with order in $\mathcal{G}_2$, then $\langle a,b,c \rangle$ is a noncollider with order in $\mathcal{G}_1$ iff $\langle a,b,c \rangle$ is a noncollider with order in $\mathcal{G}_2$.*

PROOF. This follows directly from Lemma 3.13. □

Though Proposition 3.12 appears similar to Corollary 3.14, the premise in the former assumes the two graphs are Markov equivalent, while in the latter it does not.

3.5. *Discriminating sections of a path.* It follows from Proposition 3.12 that having the same colliders with order is a necessary condition for Markov equivalence. As a step toward showing that this condition (together with the same adjacencies) is sufficient, we will show that every triple on a "minimal" $m$-connecting path has order (see Section 3.6). We first consider, in general, the relationships between different sections of a given path, where the endpoints of each section are distinguished.

Let $\boldsymbol{\pi}$ be a path with endpoints $l, r$. Let $\mathfrak{S}_{\boldsymbol{\pi}} = \{\langle x_i, b_i \rangle \mid 1 \leq i \leq m\}$ be a set of ordered pairs of vertices on $\boldsymbol{\pi}$, such that: (a) $x_i = b_j$ implies $i \neq j$, and (b) the $b_i$ are distinct and not endpoints of $\boldsymbol{\pi}$. Define a relation on the $b_i$ in $\mathfrak{S}_{\boldsymbol{\pi}}$: $b_s \prec_\pi b_t$ if $b_s$ is a nonendpoint vertex on the section $\boldsymbol{\pi}(x_t, b_t)$.

LEMMA 3.15. *With $\mathfrak{S}_{\boldsymbol{\pi}}$ and $\prec_\pi$ as defined, if $b_1 \prec_\pi \cdots \prec_\pi b_m \prec_\pi b_1$ then there exist $b_s, b_t$ such that $b_s \prec_\pi b_t \prec_\pi b_s$, $b_s$ is on $\boldsymbol{\pi}(l, x_s)$, and $b_t$ is on $\boldsymbol{\pi}(x_t, r)$.*

PROOF. It follows, from (b), that for a given $b$ there is at most one $x$ such that $\langle x, b \rangle \in \mathfrak{S}_{\boldsymbol{\pi}}$. Let $L = \{b \mid \langle x, b \rangle \in \mathfrak{S}_{\boldsymbol{\pi}}, b \text{ is on } \boldsymbol{\pi}(l, x)\}$; similarly, let $R = \{b \mid \langle x, b \rangle \in \mathfrak{S}_{\boldsymbol{\pi}}, b \text{ is on } \boldsymbol{\pi}(x, r)\}$ (see Figure 7). $R \cap L = \varnothing$ by (a) and (b). $L \neq \varnothing$, because if $b_i, b_j \in R$ and $b_i \prec_\pi b_j$ then $b_j$ is closer to $r$ than $b_i$ on $\boldsymbol{\pi}$, but if $b_1, \ldots, b_p \in R$, then $b_1 \prec_\pi \cdots \prec_\pi b_p \prec_\pi b_1$ implies that $b_1$ is closer to $r$ than $b_1$, a contradiction. Similarly, $R \neq \varnothing$. Let $b_s$ be the vertex in $L$ that is closest to $r$. Now, define $B = \{b \mid \langle x, b \rangle \in \mathfrak{S}_{\boldsymbol{\pi}}, b \prec_\pi b_s\}$; $B \neq \varnothing$, because $b_{s^*} \prec_\pi b_s$ where $s^* = s - 1 (\operatorname{mod} m)$. By definition of $b_s$, $B \subseteq R$. Let $X = \{x \mid \text{for some } b \in B, \langle x, b \rangle \in \mathfrak{S}_{\boldsymbol{\pi}}\}$. Let $x_t$ be the vertex in $X$ that is closest to $l$, and $b_t$ be a corresponding vertex in $B$, so that $\langle x_t, b_t \rangle \in \mathfrak{S}_{\boldsymbol{\pi}}$. It is sufficient to prove that $x_t$ is on $\boldsymbol{\pi}(l, b_s)$, but $x_t \neq b_s$, since then $b_s \prec_\pi b_t \prec_\pi b_s$ as required (see Figure 7). Suppose, for a contradiction, that $x_t$ is on $\boldsymbol{\pi}(b_s, r)$.

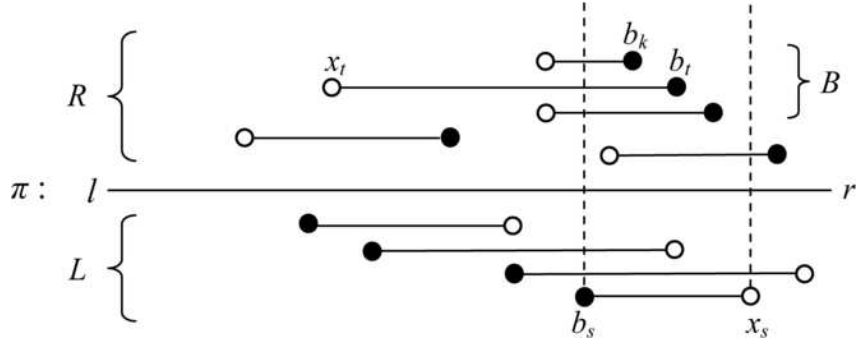

FIG. 7. *Illustration of the proof of Lemma 3.15. Lines indicate sections $\boldsymbol{\pi}(x_i, b_i)$; filled circles are $b_i$'s, open circles are $x_i$'s. Indicated are those sections for which the $b$ endpoint (filled circle) belongs to $L$, $R$ and $B$. See proof for further explanation.*

Let $b_k$ be the vertex in $B$ that is closest to $b_s$ [$b_k \neq b_s$ by (b)]. Since, by hypothesis, $x_t$ is on $\boldsymbol{\pi}(b_s, r)$, it follows by definition of $x_t$, that $x_k$ is also on $\boldsymbol{\pi}(b_s, r)$. By hypothesis, $b_{k^*} \prec_\pi b_k$ with $k^* = k - 1 \pmod{m}$. However, $b_{k^*} \notin L$, since, by definition of $b_s$, any vertex $b_i \in L$ is on $\boldsymbol{\pi}(l, b_s)$. If $b_{k^*} \in R$ then $b_{k^*} \in B$ because $x_k$ and $b_k$ are both on $\boldsymbol{\pi}(b_s, x_s)$. But then $b_k$ is not the vertex in $B$ closest to $b_s$ on $\boldsymbol{\pi}$, which is a contradiction. □

We now consider the special case of the development above, in which

(3.1) $\quad \mathfrak{S}_{\boldsymbol{\pi}} = \{\langle x_i, b_i \rangle \mid \text{ for some } y_i, \boldsymbol{\pi}(x_i, y_i) \text{ is a discriminating path for } b_i\}.$

In this context, by definition of a discriminating path, if $b_i \prec_\pi b_j$ then $b_i$ is a collider on the discriminating path $\boldsymbol{\pi}(x_j, y_j)$ for $b_j$, with $x_j, b_i, b_j$ and $y_j$ distinct vertices; $b_j$ and $y_j$ are adjacent (by the naming convention on page 10); both $b_i$ and $b_j$ are in shielded triples on $\boldsymbol{\pi}$. That $\mathfrak{S}_{\boldsymbol{\pi}}$ still satisfies (a) and (b) follows from the definition of a discriminating path together with the following.

PROPOSITION 3.16. *In a MAG $\mathcal{G}$, if $\langle a, b, c \rangle$ is a section of a path $\boldsymbol{\pi}$ between $x$ and $y$, and $a$ and $c$ are adjacent in $\mathcal{G}$, then there is at most one vertex $v$ on $\boldsymbol{\pi}$ such that either $\boldsymbol{\pi}(v, c)$ or $\boldsymbol{\pi}(v, a)$ forms a discriminating path for $b$.*

PROOF. If there is some discriminating path for $\langle a, b, c \rangle$ then $a$ is either a parent or child of $c$. In the former case, $v$ is uniquely determined as the closest vertex to $a$ on $\boldsymbol{\pi}(x, c)$ that is not a parent of $c$. The other case is symmetric: $v$ is the vertex closest to $c$ on $\boldsymbol{\pi}(a, y)$ that is not a parent of $a$. □

From here on, $\mathfrak{S}_{\boldsymbol{\pi}}$ and $\prec_\pi$ will refer to (3.1). We now prove that, as the symbol $\prec_\pi$ suggests, this relation between discriminating paths is acyclic.

COROLLARY 3.17. *On a path $\boldsymbol{\pi}$ in a MAG $\mathcal{G}$, with $\mathfrak{S}_{\boldsymbol{\pi}}$ given by (3.1), there is no sequence of distinct vertices $\langle b_1, b_2, \ldots, b_k \rangle, k > 1$, such that $b_i \prec_\pi b_{i+1}, 1 \le i < k,$ and $b_k \prec_\pi b_1$.*

This acyclic property is central to establishing that every triple on a "minimal" $m$-connecting path has order; see Lemma 3.21. Note, however, that the relation $\prec_\pi$ is not transitive in general.

PROOF OF COROLLARY 3.17. By Lemma 3.15, it is sufficient to prove that there is no pair of distinct vertices $\{b_1, b_2\}$ such that $b_1 \prec_\pi b_2$ and $b_2 \prec_\pi b_1$. For a contradiction, suppose that there is such a pair $\{b_1, b_2\}$ (see Figure 8).

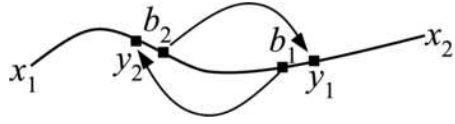

Fig. 8. *Diagram for proof of Corollary 3.17.*

By maximality, $x_1 \neq y_2$ and $x_2 \neq y_1$; otherwise, $\pi(x_1, y_1)$ or $\pi(x_2, y_2)$, respectively, would form an inducing path with nonadjacent endpoints. We now reach a contradiction because (i) $y_2$ lies on $\pi(x_1, y_1)$ and, hence, is a parent of $y_1$, but (ii) $y_1$ lies on $\pi(x_2, y_2)$ and, hence, is a parent of $y_2$. □

3.6. *Minimal m-connecting paths.* We next study the structure of "minimal" $m$-connecting paths and examine which nonconsecutive vertices on such a path may be adjacent.

DEFINITION 3.18. In a MAG, a path $\mu$, $m$-connecting $x$ and $y$ given $Z$, will be said to be *minimal* if no order preserving (proper) subsequence of the vertices on $\mu$ forms an $m$-connecting path between $x$ and $y$ given $Z$.

It is simple to see that if there is some path $m$-connecting $x$ and $y$ given $Z$, then there is a minimal path which $m$-connects $x$ and $y$ given $Z$. If $\mu = \langle v_1, \ldots, v_p \rangle$ is a path, then we will refer to any pair of vertices $(v_i, v_j)$ for which $|i - j| > 1$ as *nonconsecutive vertices on $\mu$*. As the next lemma shows, on a minimal $m$-connecting path, only certain nonconsecutive vertices may be adjacent.

LEMMA 3.19. *Let $\pi$ be a minimal m-connecting path between $a$ and $b$ given $Z$ in the MAG $\mathcal{G}$. If $i$ and $j$ are two nonconsecutive vertices on $\pi$ that are adjacent in $\mathcal{G}$ ($a = i$ or $j = b$ are possible) then exactly one of $i$ and $j$ is:* (i) *a collider on $\pi$,* (ii) *in $Z$ and* (iii) *a parent of the other vertex.*

Note that the existence of nonconsecutive vertices on a minimal $m$-connecting path implies that there are at least four vertices on the path. Lemma 3.19 is illustrated in Figure 9.

PROOF OF LEMMA 3.19. Suppose that $j$ is on $\pi(i, b)$; the other case is symmetric. Let $\eta$ be the path formed by concatenating $\pi(a, i)$ with the $\langle i, j \rangle$ edge and $\pi(j, b)$ (omit the relevant section if $i = a$ or $j = b$). Define the *status* of a vertex to be one of either an endpoint, a collider or a noncollider.

Suppose $i$ has the same status along $\eta$ as it does along $\pi$, and similarly so for $j$. Then, clearly, both $\pi$ and $\eta$ are $m$-connecting given $Z$, but $\eta$ is shorter than $\pi$, thereby violating the minimality of $\pi$. Hence, at least one of $i$ and $j$ has a status on $\boldsymbol{\eta}$ different from that on $\boldsymbol{\pi}$. Without loss of generality, suppose it is $i$; again, the other case is symmetric. $i$ is not an endpoint, because $\boldsymbol{\pi}$ and $\boldsymbol{\eta}$ have the same endpoints. It follows that *either* $i$ is a collider on $\boldsymbol{\eta}$ and $i \notin \mathrm{an}(Z)$, *or* $i$ is a noncollider on $\boldsymbol{\eta}$ and $i \in Z$.

Suppose the former, so $i \notin \mathrm{an}(Z)$, $i$ is a collider along $\boldsymbol{\eta}$, but $i$ is a noncollider along $\boldsymbol{\pi}$. Since $i$ is a collider on $\boldsymbol{\eta}$, and $\boldsymbol{\pi}(a,i) = \boldsymbol{\eta}(a,i)$, there is an arrowhead at $i$ on $\boldsymbol{\pi}(a,i)$. Then by Lemma 2.4, since $i \notin \mathrm{an}(Z)$, $\boldsymbol{\pi}(i,b)$ forms a directed path from $i$ to $b$. But $j$ is on $\boldsymbol{\pi}(i,b)$, and $i$ is a collider on $\boldsymbol{\eta}$; hence, $j\text{?}\!\!\rightarrowtail\!\! i \!\rightarrowtail\!\cdots\!\rightarrowtail\! j$ which violates Definition 2.1(a), (b).

Hence, $i \in Z$, $i$ is a noncollider along $\boldsymbol{\eta}$, but $i$ is a collider along $\boldsymbol{\pi}$. Thus, $i\text{—?}j$ in $\mathcal{G}$. Finally, the $\langle i,j \rangle$ edge cannot be undirected because $a\text{?–?}\cdots\text{?}\!\!\rightarrowtail\!\! i\!\!-\!\!\!-\!\!\!-\!\! j$ violates Definition 2.1(c); hence, $i \!\rightarrowtail\! j$. $\square$

3.7. *Discriminating paths on minimal $m$-connecting paths.* The next lemma shows that, if a triple $\langle d,b,y \rangle$ on a minimal $m$-connecting path $\boldsymbol{\pi}$ is shielded, then a subsequence of the path forms a discriminating path for $b$. Thus, in the notation of Section 3.5 on a minimal $m$-connecting path in a MAG, the following holds:

$\langle d,b,y \rangle$ a shielded triple on $\boldsymbol{\pi}$ $\implies$ there exists a nonendpoint vertex $a$

on $\boldsymbol{\pi}$ such that $a \prec_\pi b$.

LEMMA 3.20. *Let $\boldsymbol{\pi}$ be a minimal $m$-connecting path between $u$ and $v$ given $Z$ in the MAG $\mathcal{G}$. If $\langle x,b,y \rangle$ is a triple along $\boldsymbol{\pi}$ and $x$ is adjacent to $y$, then $\boldsymbol{\pi}$ contains a unique section that forms a discriminating path for $b$.*

It follows, from Lemma 3.19, that, with the possible exception of $b$, every nonendpoint vertex on the section forming a discriminating path is in $Z$.

PROOF OF LEMMA 3.20. Suppose, for a contradiction, that no such unique section exists. By Lemma 3.19, at least one of $x$ and $y$ is: (i) a collider along $\boldsymbol{\pi}$, (ii) a vertex in $Z$ and (iii) a parent of the other vertex. Without loss of generality, suppose $x$ is the vertex satisfying (i), (ii) and (iii). Since $x$ is a collider on $\pi$, we have $x\mathrel{\ast\!\!-\!\!\ast} b$ in $\mathcal{G}$. Further, since $b\mathrel{?\!\!\rightarrow} x\mathrel{\rightarrow\!\!\ast} y\mathrel{\ast\!\!-\!\!\ast} b$ in $\mathcal{G}$, by Lemma 2.2 we have $b\mathrel{?\!\!\rightarrow\!\!\ast} y$, as shown in Figure 10.

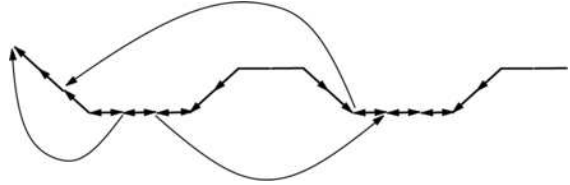

FIG. 9. *Example of a minimal $m$-connecting path (indicated by thicker edges). Here, $Z$ is the set of colliders on the path.*

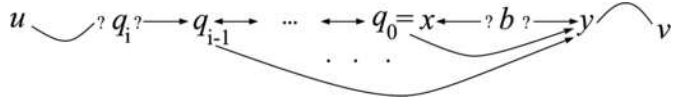

FIG. 10. *The path $\pi$ from $u$ to $v$ contains a unique section forming a discriminating path for $b$ in $\mathcal{G}$. See Lemma 3.20 for further explanation.*

Let $q_0 \equiv x$ and let $i$ be such that $q_i$ is the vertex nearest $b$ on $\pi(u,b)$ that does not satisfy at least one of the conditions (i), (ii) and (iii) satisfied by $x$. Such a vertex exists because $u$ is an endpoint and thus does not satisfy (i). Hence, $q_i$ is a vertex on $\pi(u,q_0)$ but $q_i \neq q_0$.

We now show that $q_i$ is not adjacent to $y$. Suppose otherwise. Since $q_i\mathrel{?\!\!\rightarrow\!\!\ast} q_{i-1}\mathrel{\rightarrow\!\!\ast} y$, by Lemma 2.2, we have $q_i\mathrel{?\!\!\rightarrow\!\!\ast} y$. By Lemma 3.19, (i), (ii) and (iii) are satisfied so $q_i$ is a collider on $\pi$ (hence, $q_i \neq u$), $q_i \in Z$ and $q_i\mathrel{\rightarrow\!\!\ast} y$. But this contradicts the definition of $q_i$.

Hence, $\pi(q_i, y)$ forms a discriminating path for $b$. Uniqueness follows from Proposition 3.16. □

3.8. *Triples on minimal m-connecting paths.* We now prove that, in a MAG, $\mathcal{G}$ every triple on a minimal $m$-connecting path has an order, and thus, by Proposition 3.12, is of the same type in every MAG $\mathcal{G}^*$ with $\mathcal{G}^* \sim \mathcal{G}$.

LEMMA 3.21. *If $\langle a, b, c \rangle$ is a triple on a minimal m-connecting path $\pi$ between $x$ and $y$ given $Z$ in the MAG $\mathcal{G}$, then $\langle a, b, c \rangle$ has order.*

PROOF. Suppose, for a contradiction, that $\langle a, b, c \rangle$ does not have order. Then, $a$ and $c$ are adjacent; otherwise, $\langle a, b, c \rangle$ is unshielded, and, hence, is of order 0. It follows from Lemma 3.20 that there is a unique section of $\pi$ which forms a discriminating path for $\langle a, b, c \rangle$. If every triple on this discriminating path has order, then, by definition, $\langle a, b, c \rangle$ has order. Hence, there is at least one triple which does not have order, call this $\langle a_1, b_1, c_1 \rangle$. As before, it follows that $a_1$ and $c_1$ are adjacent, and, hence, there is a unique section of $\pi$ which forms a discriminating path for $\langle a_1, b_1, c_1 \rangle$. Arguing in this way, we can construct an infinite sequence of shielded triples on $\pi$, $\langle a_i, b_i, c_i \rangle$ ($i \in \mathbb{N}$), none of which have order and such that

$$\cdots \prec_\pi b_i \prec_\pi \cdots \prec_\pi b_1 \prec_\pi b.$$

However, by Corollary 3.17 all of the $b_i$'s are distinct, which is a contradiction since $\pi$ is finite. Thus, every triple on $\pi$ has an order. □

Note that this argument shows that every triple on a minimal $m$-connecting path $\pi$ has *some* order and also that this order is bounded by the number of vertices on $\boldsymbol{\pi}$; see page 28. Though we will at no stage need to do so, note that to determine *which* order a given triple on $\boldsymbol{\pi}$ has, it might be necessary to consider other discriminating paths for the given triple, not merely those which are sections of $\boldsymbol{\pi}$.

COROLLARY 3.22. *Suppose that $\mathcal{G}_1$ and $\mathcal{G}_2$ are MAGs with the same adjacencies and the same colliders with order. If $\boldsymbol{\pi}$ is a minimal m-connecting path between $x$ and $y$ given $Z$ in $\mathcal{G}_1$, then $\langle a, b, c \rangle$ is a collider [noncollider] on $\boldsymbol{\pi}$ in $\mathcal{G}_1$ if and only if $\langle a, b, c \rangle$ is a collider [noncollider] on the corresponding path $\boldsymbol{\pi}^*$ in $\mathcal{G}_2$.*

PROOF. This follows directly from Corollary 3.14 and Lemma 3.21. □

3.9. *Directed paths from colliders to vertices in $Z$.* In this section, we establish that if there is an $m$-connecting path $\widetilde{\boldsymbol{\pi}}$ between $x$ and $y$ given $Z$ in $\mathcal{G}$, then we can always find a path $\boldsymbol{\pi}$ $m$-connecting $x$ and $y$ given $Z$ in $\mathcal{G}$ such that, if $c$ is a collider on $\boldsymbol{\pi}$, then $c$ is an ancestor of a vertex in $Z$ in any graph $\mathcal{G}^*$ which contains the same adjacencies and the same colliders with order as $\mathcal{G}$.

Let $|\boldsymbol{\pi}|$ be the *length* of a path (i.e., the number of edges on $\boldsymbol{\pi}$). Let $\mathfrak{D}(b, Z)$ be the set of directed paths from $b$ to some vertex in $Z$. $\widetilde{\boldsymbol{\delta}} \in \mathfrak{D}(b, Z)$ is said to be *a minimal directed path with respect to $Z$* if $|\widetilde{\boldsymbol{\delta}}| = \min_{\boldsymbol{\delta} \in \mathfrak{D}(b,Z)} |\boldsymbol{\delta}|$. Let

$$\phi(b, Z) = \begin{cases} 0, & \text{if } b \in Z, \\ \min_{\boldsymbol{\delta} \in \mathfrak{D}(b,Z)} |\boldsymbol{\delta}|, & \text{if } b \in \mathrm{an}(Z) \setminus Z. \end{cases}$$

If $\boldsymbol{\pi}$ $m$-connects given $Z$, then let

$$\phi(\boldsymbol{\pi}, Z) = \sum_{b \text{ a collider on } \boldsymbol{\pi}} \phi(b, Z).$$

We now construct an ordering on the set of paths $m$-connecting given $Z$:

$$\boldsymbol{\pi}_1 \ll_Z \boldsymbol{\pi}_2 \iff |\boldsymbol{\pi}_1| < |\boldsymbol{\pi}_2| \text{ or}$$
$$|\boldsymbol{\pi}_1| = |\boldsymbol{\pi}_2| \text{ and } \phi(\boldsymbol{\pi}_1, Z) < \phi(\boldsymbol{\pi}_2, Z).$$

DEFINITION 3.23. In an ancestral graph, an $m$-connecting path $\boldsymbol{\pi}$ between $x$ and $y$ given $Z$ is said to be a *closest m-connecting path to $Z$* if there is no other path $\boldsymbol{\pi}^*$ $m$-connecting $x$ and $y$ given $Z$ such that $\boldsymbol{\pi}^* \ll_Z \boldsymbol{\pi}$.

PROPOSITION 3.24. *In an ancestral graph, if there is an m-connecting path $\widetilde{\boldsymbol{\pi}}$ between $x$ and $y$ given $Z$, then there is an m-connecting path $\boldsymbol{\pi}$ that is closest to $Z$. Every such path is also a minimal m-connecting path given $Z$.*

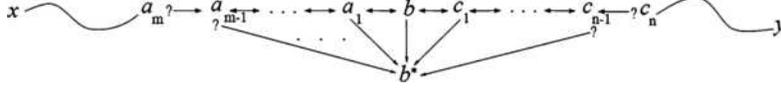

FIG. 11. *Diagram for the proof of Lemma 3.25. Either $\langle x,\ldots,a_{m-1},b^*,c_{n-1},\ldots,y\rangle$ is an $m$-connecting path closer to $Z$, or at least one of the noncolliders $\langle a_1,b,b^*\rangle$ and $\langle c_1,b,b^*\rangle$ has order. (Note that $a_0 = b = c_0$ and $m,n > 0$ by construction.)*

PROOF. Existence of a closest path $\pi$ is immediate since $\ll_Z$ is an ordering on the finite and nonempty (by hypothesis) set of paths $m$-connecting $x$ and $y$ given $Z$. Minimality follows, because if there were an $m$-connecting path $\pi^*$ formed from an order preserving (proper) subsequence of the vertices on $\pi$ then $|\pi^*| < |\pi|$, so $\pi^* \ll_Z \pi$, which is a contradiction. □

LEMMA 3.25. *If, in a MAG $\mathcal{G}$: $\pi = \langle x,\ldots,y\rangle$ is a closest $m$-connecting path to $Z$; $\langle a_1,b,c_1\rangle$ is a collider on $\pi$; and $\boldsymbol{\delta} = \langle b,b^*,\ldots,z\rangle$ is a minimal directed path with respect to $Z$ from $b$ to some $z \in Z$; then at least one of the noncolliders $a_1\mathord{?}\mathord{\rightarrow}b\mathord{\rightarrow}b^*$ or $b^*\mathord{\leftarrow}b\mathord{\leftarrow}\mathord{?}c_1$ has order in $\mathcal{G}$.*

PROOF. By Proposition 3.24, $\pi$ is a minimal $m$-connecting path between $x$ and $y$ given $Z$. Now, suppose for a contradiction that neither triple $a_1\mathord{?}\mathord{\rightarrow}b\mathord{\rightarrow}b^*$ nor $b^*\mathord{\leftarrow}b\mathord{\leftarrow}\mathord{?}c_1$ has order. Then, $a_1$ is adjacent to $b^*$ and by Lemma 2.2 we have $a_1\mathord{?}\mathord{\rightarrow}b^*$. Similarly, $c_1\mathord{?}\mathord{\rightarrow}b^*$.

Define $a_0 = b$, and let $a_m$ be the vertex along $\pi(x,b)$ that is furthest from $b$ such that for all $k$, $0 \le k < m$: (i) $a_k$ is a collider on $\pi$, and (ii) $a_k\mathord{\rightarrow}b^*$ (see Figure 11). Such a vertex $a_m$ exists because $a_0 = b$ satisfies the conditions for $a_k$; note that $m > 0$. Then, the following hold:

(1) $a_m$ is adjacent to $b^*$. Otherwise for $m = 1$, $\langle a_1,b,b^*\rangle$ is unshielded; or for $m > 1$, $\langle a_m,\ldots,a_1,b,b^*\rangle$ forms a discriminating path with order for $\langle a_1,b,b^*\rangle$ (by Lemma 3.21). In either case, $\langle a_1,b,b^*\rangle$ would have order which is a contradiction.
(2) Since $a_m\mathord{?}\mathord{\rightarrow}a_{m-1}\mathord{\rightarrow}b^*$, by Lemma 2.2, we have that $a_m\mathord{?}\mathord{\rightarrow}b^*$ is in $\mathcal{G}$.
(3) If $a_m \ne x$, then triples $\langle a_{m+1},a_m,a_{m-1}\rangle$ and $\langle a_{m+1},a_m,b^*\rangle$ are of the same type (collider/noncollider) where $a_{m+1}$ is the predecessor of $a_m$ along $\pi(x,a_m)$: if $a_m\mathord{\leftrightarrow}b^*$ then, since $a_m\mathord{?}\mathord{\rightarrow}a_{m-1}\mathord{\rightarrow}b^*$, by Lemma 2.2 we have that $a_m\mathord{\leftrightarrow}a_{m-1}$. If $a_m\mathord{\rightarrow}b^*$, then by the definition of $a_m$, triple $\langle a_{m+1},a_m,a_{m-1}\rangle$ is not a collider.

Define $c_0 = b$. Let $c_n$ be the vertex along $\pi(b,y)$ that is furthest from $b$ such that, for all $j$, $0 \le j < n$: (i) $c_j$ is a collider on $\pi$, and (ii) $c_j\mathord{\rightarrow}b^*$. By symmetric arguments to (1), (2) and (3), we may show that $c_n\mathord{?}\mathord{\rightarrow}b^*$, and either $c_n = y$ or the triples $\langle c_{n-1},c_n,c_{n+1}\rangle$ and $\langle b^*,c_n,c_{n+1}\rangle$ are of the same type, where $c_{n+1}$ is the successor of $c_n$ on the path $\pi(b,y)$ (see Figure 11).

Let $\boldsymbol{\eta}$ be the path formed by concatenating the section $\boldsymbol{\pi}(x, a_m)$ to $a_m {?}{\rightarrowtail} b^*$ $\leftarrowtail {?} c_n$ and $\boldsymbol{\pi}(c_n, y)$ (if $x = a_m$, or $c_n = y$ then omit the relevant sections). $\boldsymbol{\eta}$ forms an $m$-connecting path given $Z$ because $a_m$ and $c_n$ have the same status on $\boldsymbol{\eta}$ as they have on $\boldsymbol{\pi}$, and $b^*$ is an ancestor of $Z$. However, since $|\boldsymbol{\eta}| \leq |\boldsymbol{\pi}|$ and $\phi(\boldsymbol{\eta}, Z) < \phi(\boldsymbol{\pi}, Z)$, $\boldsymbol{\eta} \ll_Z \boldsymbol{\pi}$, which is a contradiction. □

LEMMA 3.26. *In a MAG $\mathcal{G}$ if $\boldsymbol{\delta}$ is a directed path from $v$ to $z \in Z$ and $\boldsymbol{\delta}$ is minimal with respect to $Z$, then every noncollider on $\boldsymbol{\delta}$ is unshielded (= order 0).*

PROOF. Suppose that $\langle a, b, c \rangle$ is a noncollider on $\boldsymbol{\delta}$ and $a {\rightarrowtail} b {\rightarrowtail} c$. If $a$ and $c$ are adjacent then, by Definition 2.1(a), (b), we have $a {\rightarrowtail} c$, which contradicts the minimality of $\boldsymbol{\delta}$, with respect to $Z$. □

Though not needed, in fact no nonconsecutive vertices on $\boldsymbol{\delta}$ are adjacent.

COROLLARY 3.27. *Let $\mathcal{G}_1$, $\mathcal{G}_2$ be MAGs with the same adjacencies, and the same colliders with order. If in $\mathcal{G}_1$: $\boldsymbol{\pi}$ $m$-connects $x$ and $y$ given $Z$; $\boldsymbol{\pi}$ is a closest path to $Z$; $\langle a, b, c \rangle$ is a collider on $\boldsymbol{\pi}$; and $\boldsymbol{\delta}$ forms a directed path from $b$ to a vertex $z \in Z$ that is minimal with respect to $Z$; then the corresponding path $\boldsymbol{\delta}^*$ is a directed path in $\mathcal{G}_2$.*

PROOF. By Proposition 3.24, $\boldsymbol{\pi}$ is a minimal $m$-connecting path. Let $\boldsymbol{\delta} = \langle b = d_0, \ldots, d_n = z \rangle$. The proof is by induction on the edges $\langle d_i, d_{i+1} \rangle$ of $\boldsymbol{\delta}^*$.

*Base case* ($i = 0$). Since $\langle a, b, c \rangle$ is a collider on $\boldsymbol{\pi}$ in $\mathcal{G}_1$, and $\boldsymbol{\pi}$ is minimal, by Corollary 3.22, $\langle a, b, c \rangle$ is also a collider in $\mathcal{G}_2$. By Lemma 3.25 at least one of the noncolliders, $\langle a, b, d_1 \rangle$, $\langle c, b, d_1 \rangle$ has order in $\mathcal{G}_1$, and by Corollary 3.14 is also a noncollider in $\mathcal{G}_2$. It follows by Definition 2.1(c) that $b {\rightarrowtail} d_1$ in $\mathcal{G}_2$ as required (so in fact $\langle a, b, d_1 \rangle$ and $\langle d_1, b, c \rangle$ are both noncolliders in $\mathcal{G}_2$).

*Inductive case* ($1 \leq i < n$). Assume that the section $\boldsymbol{\delta}^*(b, d_i)$ forms a directed path from $b$ to $d_i$ in $\mathcal{G}_2$. By Lemma 3.26 the noncollider $\langle d_{i-1}, d_i, d_{i+1} \rangle$ has order, and hence is a noncollider in $\mathcal{G}_2$. By the induction hypothesis we have $d_{i-1} {\rightarrowtail} d_i$ in $\mathcal{G}_2$; hence, by Definition 2.1(c), $d_i {\rightarrowtail} d_{i+1}$ in $\mathcal{G}_2$, as required. □

3.10. *Characterization of Markov equivalence.* We now prove the main result of this paper, Theorem 3.7.

PROOF OF THEOREM 3.7. (*if*) Since $\mathcal{G}_1$ and $\mathcal{G}_2$ have the same adjacencies and colliders with order, by Corollary 3.14, $\mathcal{G}_1$ and $\mathcal{G}_2$ also have the same noncolliders with order. By definition, $X$ is $m$-separated from $Y$ given $Z$ if and only if for all $x \in X$, $y \in Y$, $x$ is $m$-separated from $y$ given $Z$. Thus, it is sufficient to show that $x$ and $y$ are $m$-connected given $Z$ in $\mathcal{G}_1$ if and only if $x$ and $y$ are $m$-connected given $Z$ in $\mathcal{G}_2$. If $x$ and $y$ are $m$-connected given $Z$ in $\mathcal{G}_1$, then, by Proposition 3.24, there exists a path $\boldsymbol{\pi}$ which $m$-connects $x$ and $y$ given $Z$, is minimal and is closest to $Z$ in $\mathcal{G}_1$. By Corollary 3.22 every triple on $\boldsymbol{\pi}$ is of the same type on the corresponding path $\boldsymbol{\pi}^*$ in $\mathcal{G}_2$. Hence, every noncollider on $\boldsymbol{\pi}^*$ is not in $Z$. Since $\boldsymbol{\pi}$ is $m$-connecting, every collider $b$ on $\boldsymbol{\pi}$ is an ancestor of $Z$; hence, if $b \notin Z$ then there exists a directed path $\boldsymbol{\delta}_b$ from $b$ to some vertex $z_b \in Z$ that is minimal with respect to $Z$. By Corollary 3.27, the corresponding path $\boldsymbol{\delta}_b^*$ forms a directed path from $b$ to $z_b$ in $\mathcal{G}_2$. Thus, every collider on $\boldsymbol{\pi}^*$ is an ancestor of $Z$ in $\mathcal{G}_2$ and $\boldsymbol{\pi}^*$ $m$-connects $x$ and $y$ given $Z$ in $\mathcal{G}_2$. Likewise, it is easy to see (by symmetry) that an $m$-connecting path in $\mathcal{G}_2$ implies that there is an $m$-connecting path in $\mathcal{G}_1$. Thus, $\mathcal{G}_1$ and $\mathcal{G}_2$ are Markov equivalent.

(*only if*) Conversely, if $\mathcal{G}_1$ and $\mathcal{G}_2$ are Markov equivalent, then, by Proposition 3.6, they have the same adjacencies, and, by Proposition 3.12, they have the same colliders with order. □

COROLLARY 3.28. *Two ancestral graphs $\mathcal{G}_1$ and $\mathcal{G}_2$ are Markov equivalent iff the corresponding unique MAGs $\bar{\mathcal{G}}_1$ and $\bar{\mathcal{G}}_2$ of which $\mathcal{G}_1$ and $\mathcal{G}_2$ are, respectively, subgraphs and to which they are Markov equivalent, satisfy the conditions given in Theorem 3.7.*

**4. Related work and computational complexity.** Two prior characterizations of Markov equivalence for MAGs have been given in the literature.

THEOREM 4.1 [19]. *Two MAGs $\mathcal{G}_1$ and $\mathcal{G}_2$ are Markov equivalent if and only if:*

  (i) *$\mathcal{G}_1$ and $\mathcal{G}_2$ have the same adjacencies;*
  (ii) *$\mathcal{G}_1$ and $\mathcal{G}_2$ have the same unshielded colliders; and*
  (iii) *if $\boldsymbol{\pi}$ forms a discriminating path for $b$ in $\mathcal{G}_1$ and $\mathcal{G}_2$, then $b$ is a collider on $\boldsymbol{\pi}$ in $\mathcal{G}_1$ if and only if it is a collider on $\boldsymbol{\pi}$ in $\mathcal{G}_2$.*

More recently, [24] gave the following elegant characterization.

THEOREM 4.2 [24]. *Two MAGs $\mathcal{G}_1$ and $\mathcal{G}_2$ are Markov equivalent if and only if $\mathcal{G}_1$ and $\mathcal{G}_2$ have the same minimal collider paths.*

Here, a collider path $\boldsymbol{\nu} = \langle v_1, \ldots, v_n \rangle$ is *minimal* if there is no order preserving subsequence $\langle v_1 = v_{i_1}, \ldots, v_{i_k} = v_n \rangle$, which forms a collider path (single edges are trivially minimal collider paths).

However, neither of these characterizations lead to a polynomial time algorithm. Clause (iii) in Theorem 4.1 requires us to verify that, if there is a discriminating path in both $\mathcal{G}_1$ and $\mathcal{G}_2$, then the triple discriminated is a collider or noncollider in both. Thus, in principle, we need to find every discriminating path for a given triple; otherwise, it is possible that, although a triple is discriminated by *some* path in $\mathcal{G}_1$ and *some* path in $\mathcal{G}_2$, in fact there is no discriminating path that is common to both graphs. Since the number of such paths may grow at super-polynomial rate, finding them all would not be feasible in polynomial-time. (Reference [19] outlined a method for checking Markov equivalence using the conditions of Theorem 3.7, rather than Theorem 4.1, though the paper only proves the latter result. The computational complexity claim in that paper was also incorrect.)

Similarly, it is not hard to show that the number of minimal collider paths in a graph may grow super-polynomially with the number of vertices so the conditions in Theorem 4.2 cannot, in general, be verified in polynomial time (see supplementary material [1]).

In the Appendix, we provide an algorithm that verifies the conditions in Theorem 3.7 in $O(ne^4)$ calculations, where the graphs have $n$ vertices, and $e$ edges. For a general, not necessarily maximal, ancestral graph $\mathcal{G}$ the unique MAG $\bar{\mathcal{G}}$ of which it is a subgraph and to which it is Markov equivalent may be found in $O(n^5)$ time; thus, the conditions in Corollary 3.28 may also be checked in polynomial time.

4.1. *Summary graphs and MC graphs.* *Summary graphs*, described in Cox and Wermuth [5], represent another approach to representing the independence structure of DAGs under marginalizing and conditioning. For a given summary graph $\mathcal{H}$, it is always possible to construct a DAG $\mathcal{D}(H)$ with additional variables such that the DAG is Markov equivalent to $\mathcal{H}$ after marginalizing and conditioning. Consequently, it is always possible to transform a summary graph into an ancestral graph via the graphical transformation mentioned in Section 2.6. Hence, via this transformation, the results in this paper also provide an algorithm for determining the Markov equivalence of two summary graphs. We note that in general it may not be possible to recover the summary graph from the corresponding ancestral graph (see [15], Section 9).

Koster introduced another class of graphs, called *MC-graphs*, together with an operation of marginalizing and conditioning (see [10, 11]). For MC-graphs it is not always the case that there exists some DAG which is Markov equivalent to the MC-graph under marginalizing and conditioning. However, for the subclass of MC-graphs which are Markov equivalent to DAGs with additional variables under marginalizing and conditioning, we may again apply the results of this paper to establish Markov equivalence.

TABLE A.1
*The algorithm* **Reachable**$(\mathbb{D}, \boldsymbol{w})$

| **Inputs**: | a directed graph $\mathbb{D}(\mathbb{V}, \mathbb{E})$; an element $\boldsymbol{w} \in \mathbb{V}$ |
|---|---|
| **Output**: | a set $\mathbb{S}$ of elements connected to $\boldsymbol{w}$ in $\mathbb{D}$ |
| 1 | $\mathbb{S}_0 = \varnothing; \mathbb{S}_1 = \{\boldsymbol{w}\}; p = 1;$ |
| 2 | `repeat` |
| 3 | $\mathbb{S}_{p+1} = \mathbb{S}_p \cup \{\boldsymbol{w}_2 | \boldsymbol{w}_1 \in \mathbb{S}_p \setminus \mathbb{S}_{p-1} \text{ and } \langle \boldsymbol{w}_1, \boldsymbol{w}_2 \rangle \in \mathbb{E}\};$ |
| 4 | $p = p + 1;$ |
| 5 | `until` $\mathbb{S}_p = \mathbb{S}_{p-1};$ |
| 6 | `return` $\mathbb{S} = \mathbb{S}_p.$ |

## APPENDIX

We introduce the following notation:

$$\mathfrak{Adj}(\mathcal{G}) = \{\langle x, y \rangle \mid x \text{ and } y \text{ are adjacent in } \mathcal{G}\},$$

$$\mathfrak{Col}(\mathcal{G}) = \{\langle x, y, z \rangle \mid x\text{?}{\rightarrowtail}y{\leftarrowtail}\text{?}z \text{ in } \mathcal{G}\},$$

$$\mathfrak{OCol}(\mathcal{G}) = \{\langle x, y, z \rangle \mid \langle x, y, z \rangle \in \mathfrak{Col}(\mathcal{G}) \text{ and } \langle x, y, z \rangle \text{ has order}\},$$

$$\mathfrak{ICol}(\mathcal{G}) = \bigcap_{\mathcal{G}^* \sim \mathcal{G}} \mathfrak{Col}(\mathcal{G}^*),$$

which are, respectively, the set of adjacencies, colliders, colliders with order in $\mathcal{G}$ and colliders common to all graphs in the Markov equivalence class containing $\mathcal{G}$. In general, we have $\mathfrak{OCol}(\mathcal{G}) \subseteq \mathfrak{ICol}(\mathcal{G}) \subseteq \mathfrak{Col}(\mathcal{G})$.

TABLE A.2
*The algorithm* **Triples**$(\mathcal{G})$

| **Input**: | a maximal ancestral graph $\mathcal{G}$ |
|---|---|
| **Output**: | a set of triples $\mathbb{T}$ such that $\mathfrak{OCol}(\mathcal{G}) \subseteq \mathbb{T} \subseteq \mathfrak{ICol}(\mathcal{G})$ |
| 1 | $\mathbb{T}_0 = \{\langle a, b, c \rangle | \langle a, b, c \rangle \in \mathfrak{Col}(\mathcal{G}), (a, c) \notin \mathfrak{Adj}(\mathcal{G})\};$ |
| 2 | $k = 0;$ |
| 3 | `repeat` |
| 4 | $k = k + 1; \mathbb{T}_k = \mathbb{T}_{k-1};$ |
| 5 | `for each` $\langle a, b, c \rangle \in \mathfrak{Col}(\mathcal{G}) \setminus \mathbb{T}_{k-1}$ `with` $a \in \text{sp}_{\mathcal{G}}(b) \cap \text{pa}_{\mathcal{G}}(c):$ |
| 6 | $\mathbb{V} = \{\langle t, u \rangle | t, u \in \text{pa}(c), t {\leftarrowtail}{\rightarrowtail} u \text{ in } \mathcal{G}\} \cup \{\langle b, a \rangle\};$ |
| 7 | $\mathbb{E} = \{\langle \langle t, u \rangle, \langle u, v \rangle \rangle | \langle t, u, v \rangle \in \mathbb{T}_{k-1}, \langle t, u \rangle, \langle u, v \rangle \in \mathbb{V}\};$ |
| 8 | $\mathbb{S} = \textbf{Reachable}((\mathbb{V}, \mathbb{E}), \langle b, a \rangle);$ |
| 9 | $X = \{x \mid \exists y, z, \langle z, y, x \rangle \in \mathbb{T}_{k-1}, \langle z, y \rangle \in \mathbb{S}\};$ |
| 10 | `if` $X \setminus \{v \mid (v, c) \in \mathfrak{Adj}(\mathcal{G})\} \neq \varnothing$ |
| | `then` $\mathbb{T}_k = \mathbb{T}_k \cup \{\langle a, b, c \rangle, \langle c, b, a \rangle\};$ |
| 11 | `until` $\mathbb{T}_k = \mathbb{T}_{k-1};$ |
| 12 | `return` $\mathbb{T} = \mathbb{T}_k.$ |

TABLE A.3
*The algorithm* **Equivalent**($\mathcal{G}_1$, $\mathcal{G}_2$)

| | |
|---|---|
| **Inputs**: | two maximal ancestral graphs $\mathcal{G}_1$ and $\mathcal{G}_2$ |
| **Output**: | a Boolean variable indicating whether $\mathfrak{I}_m(\mathcal{G}_1) = \mathfrak{I}_m(\mathcal{G}_2)$ |
| 1 | if $\mathfrak{Adj}(\mathcal{G}_1) \neq \mathfrak{Adj}(\mathcal{G}_2)$ return FALSE; |
| 2 | if $\mathbf{Triples}(\mathcal{G}_1) \setminus \mathfrak{Col}(\mathcal{G}_2) \neq \varnothing$ return FALSE; |
| 3 | if $\mathbf{Triples}(\mathcal{G}_2) \setminus \mathfrak{Col}(\mathcal{G}_1) \neq \varnothing$ return FALSE; |
| 4 | return TRUE. |

The equivalence algorithm is described in Tables A.1–A.3. The main procedure, **Triples**($\mathcal{G}$), identifies a superset of the colliders with order as follows. A discriminating path $\pi = \langle x, z, \ldots, a, b, c \rangle$ for the collider $\langle a, b, c \rangle$ (where $z$ may equal $a$) that is in $\mathbb{T}_k$ divides naturally into three parts. First, there is a collider $\langle a, b, c \rangle$, which is not in $\mathbb{T}_{k-1}$. Second, there is a collider path $\gamma = \langle z \equiv v_1 \leftrightarrow \cdots \leftrightarrow v_j \equiv b \rangle$, where $v_1, \ldots, v_j \in \mathrm{pa}(c)$, and the triples $\langle v_{i-1}, v_i, v_{i+1} \rangle \in \mathbb{T}_{k-1}$. The third part is an edge $x ? \rightarrow z$, where $x$ is not adjacent to $c$ and for some $y$, $x ? \rightarrow z \leftrightarrow y \in \mathbb{T}_{k-1}$, and $z \leftrightarrow y$ is on the path $\gamma$. Line 5 of **Triples**($\mathcal{G}$) locates candidate triples $\langle a, b, c \rangle$. Steps 6, 7, and 8 search for collider paths $\gamma$. Note that "vertices" ($\mathbb{V}$) and "edges" ($\mathbb{E}$) in $\mathbb{D}$ correspond to, respectively, edges and colliders in $\mathcal{G}$. Finally, lines 9 and 10 search for a vertex satisfying the conditions on $x$. For further insight into the operation of the algorithm, we refer the reader to the proof of correctness.

PROPOSITION A.1. *The algorithm* **Triples**($\mathcal{G}$) *returns a set* $\mathbb{T}$ *satisfying* (a) $\mathfrak{OCol}(\mathcal{G}) \subseteq \mathbb{T}$ *and* (b) $\mathbb{T} \subseteq \mathfrak{ICol}(\mathcal{G})$.

PROOF OF (a). The proof is by induction on the order of the collider. By construction, $\mathbb{T}_0$ is the set of unshielded colliders in $\mathcal{G}$, which is the set of colliders of order 0. Our induction hypothesis is that all colliders with order less than $k > 0$ are contained in $\mathbb{T}_{k-1}$, at line 11. If $\langle a, b, c \rangle$ is a collider with order $k$, then either $a \rightarrow c$ or $c \rightarrow a$. Suppose the former. Then, there exists a discriminating path $\langle x = q_0, q_1, \ldots, q_p = a, q_{p+1} = b, c \rangle$ on which $\langle q_{j-1}, q_j, q_{j+1} \rangle$ ($1 \leq j \leq p$) are colliders of order less than $k$. By definition of a discriminating path, $\langle q_{p-1}, a, b \rangle$ is a collider, as is $\langle a, b, c \rangle$, so $a \in \mathrm{sp}(b)$. Thus, $\langle a, b, c \rangle$ satisfies the conditions at line 5. In addition, for $1 \leq j \leq p-1$, $q_j, q_{j+1} \in \mathrm{pa}(c)$ and $q_j \leftrightarrow q_{j+1}$, so $\langle q_{j+1}, q_j \rangle \in \mathbb{V}$. In addition, $\langle q_{p+1}, q_p \rangle = \langle b, a \rangle \in \mathbb{V}$ by construction. Since for $1 < j \leq p$, $\langle q_{j-1}, q_j, q_{j+1} \rangle$ is a collider of order less than $k$, it follows by the induction hypothesis that $\langle q_{j-1}, q_j, q_{j+1} \rangle \in \mathbb{T}_{k-1}$. Thus, $\langle \langle q_{j+1}, q_j \rangle, \langle q_j, q_{j-1} \rangle \rangle \in \mathbb{E}$. Consequently, $\langle q_2, q_1 \rangle \in \mathbb{S}$ at line 8, since the sequence $\langle \langle b = q_{p+1}, a = q_p \rangle, \langle q_p, q_{p-1} \rangle, \ldots, \langle q_2, q_1 \rangle \rangle$ is found (recursively) by calls to **Reachable**. Since $\langle q_2, q_1, x = q_0 \rangle \in \mathbb{T}_{k-1}$, it follows that $x \in X$. Finally, by definition of a discriminating path, $x$ is not adjacent to $c$. Thus, the condition in the `if` clause at line 10 holds, so if $\langle a,b,c\rangle \notin \mathbb{T}_{k-1}$ then it is added to $\mathbb{T}_k$.

PROOF OF (b). The proof is by induction on $k$ in the algorithm. We show that $\mathbb{T}_k \subseteq \mathfrak{ICol}(\mathcal{G})$. When $k=0$, $\mathbb{T}_0$ is the set of unshielded colliders, so the result follows from Proposition 3.6. For $k>0$ our induction hypothesis is that $\mathbb{T}_{k-1} \subseteq \mathfrak{ICol}(\mathcal{G})$. If $\langle a,b,c\rangle \in \mathbb{T}_k \setminus \mathbb{T}_{k-1}$, then either $\langle a,b,c\rangle$ or $\langle c,b,a\rangle$ (but not both) satisfies the condition at line 5. Suppose the former; the other case is symmetric. There exists a triple $\langle x,y,z\rangle \in \mathbb{T}_{k-1}$, with $\langle y,z\rangle \in \mathbb{S}$, and $x$ not adjacent to $c$. Since $\langle y,z\rangle \in \mathbb{S}$, and $x \in X$, there exists a sequence of edges, $\mathfrak{s} \equiv \langle\langle b,a\rangle,\ldots,\langle z,y\rangle,\langle y,x\rangle\rangle$ such that each consecutive pair of edges in $\mathfrak{s}$ forms a collider in $\mathbb{T}_{k-1}$, all vertices other than $b$ and $x$ are parents of $c$, and all edges other than possibly $\langle y,x\rangle$ are bi-directed in $\mathcal{G}$. Note that it follows from the inductive hypothesis that all of the colliders formed by successive pairs of edges in $\mathfrak{s}$ are present in any graph $\mathcal{G}^*$ Markov equivalent to $\mathcal{G}$. We have thus established that, with the possible exception of the first and last edge in the sequence, all these edges are bi-directed in every graph in the Markov equivalence class. However, the sequence of edges in $\mathfrak{s}$ may not form a path because the associated sequence of vertices may contain repeats. Removing loops leads to a unique path $\boldsymbol{\pi}$ with endpoints $b$ and $x$. By construction, $b$ and $x$ only occur in the edges $\langle b,a\rangle$ and $\langle y,x\rangle$, respectively (since $b,x$ are not parents of $c$, while all other vertices in the sequence are); consequently, these edges are on $\boldsymbol{\pi}$. Hence, $\boldsymbol{\pi}$ forms a collider path from $x$ to $b$, and all of the colliders on this path are present in every graph in the Markov equivalence class. By Lemma 3.10, $\boldsymbol{\pi}$ forms a discriminating path in every graph Markov equivalent to $\mathcal{G}$. Thus, by Lemma 3.9, $\langle a,b,c\rangle \in \mathfrak{ICol}(\mathcal{G})$ as required. □

Our proof establishes that all triples in **Triples**$(\mathcal{G})$ are colliders present in every graph in the Markov equivalence class containing $\mathcal{G}$, which might include some colliders that do not have order. If we were able to identify any triples in **Triples**$(\mathcal{G}) \setminus \mathfrak{OCol}(\mathcal{G})$ without increasing the complexity of the algorithm, then the algorithm could be made more efficient since it is redundant to check for the presence of such colliders in the other graph. However, we know of no examples where **Triples**$(\mathcal{G}) \setminus \mathfrak{OCol}(\mathcal{G}) \neq \varnothing$.

PROPOSITION A.2. *The algorithm* **Equivalent**$(\mathcal{G}_1,\mathcal{G}_2)$ *returns* TRUE *iff* $\mathcal{G}_1$ *and* $\mathcal{G}_2$ *are Markov equivalent.*

PROOF. "if" follows from Propositions 3.6 and A.1(b). "only if" follows from Proposition A.1(a), Lemma 3.13 and Theorem 3.7. □

**Reachable**($\mathbb{D}, w$) runs in time $O(\tilde{e})$ where $\tilde{e}$ is the number of edges in $\mathbb{D}$. The graph $\mathbb{D}$ may be represented as a list of adjacencies for each vertex, with each edge $\langle \boldsymbol{w}_1, \boldsymbol{w}_2 \rangle$ being considered at most once at line 3.

Now, consider the complexity of **Triples**($\mathcal{G}$). Let $n$ and $e$ denote, respectively, the number of vertices and edges in $\mathcal{G}$. Any triple appearing on a minimal $m$-connecting path $\boldsymbol{\pi}$ has order at most $n-3$: $\boldsymbol{\pi}$ contains at most $n$ vertices; hence, at most $n-2$ triples; all of the other discriminating paths involved are sections of $\boldsymbol{\pi}$; and unshielded triples (of which there is at least one) are of order 0. Thus, it is always sufficient for Markov equivalence to check that two graphs have triples of order less than $n$. Hence, the outer loop, at line 4, in **Triples**($\mathcal{G}$) is of complexity $O(n)$. The number of colliders in $\mathcal{G}$ is of $O(e^2)$; hence, the loop at line 5 is executed $O(e^2)$ times (for each $k$). Since $\mathbb{E}$ is of size $O(e^2)$, lines 6 to 8 are also of complexity $O(e^2)$. Finally, line 9 is $O(e^2)$ [since $\mathbb{T}_{k-1}$ is of size $O(e^2)$] and line 10 is $O(e)$. Thus, the overall complexity is $O(ne^4)$.

**Acknowledgments.** We thank three anonymous referees whose comments led to great improvements in this paper, particularly in Sections 2.4 and 3.9.

R. A. ALI
DEPARTMENT OF MATHEMATICS AND STATISTICS
UNIVERSITY OF GUELPH
GUELPH, ONTARIO N1G 2W1
CANADA
E-MAIL: aali@uoguelph.ca

T. S. RICHARDSON
DEPARTMENT OF STATISTICS
UNIVERSITY OF WASHINGTON
SEATTLE, WASHINGTON 98195-4322
USA
E-MAIL: tsr@stat.washington.edu

P. SPIRTES
DEPARTMENT OF PHILOSOPHY
CARNEGIE-MELLON UNIVERSITY
PITTSBURGH, PENNSYLVANIA 15213
USA
E-MAIL: ps7z@andrew.cmu.edu